\documentclass[11pt]{amsart}

\usepackage{amssymb,amsfonts,amsmath,amscd} 
\usepackage{amsthm}
\usepackage{stmaryrd}
\usepackage{graphicx,graphics}
\usepackage{latexsym}
\usepackage{color} 
\usepackage[all]{xy}
\usepackage{pstricks}
\usepackage{pst-node}
\usepackage{verbatim}
\usepackage{hyperref} 

\newtheorem{thm}{Theorem}[section]
\newtheorem{prop}[thm]{Proposition}

\newtheorem{cla}[thm]{Claim}

\newtheorem{prop-defn}[thm]{Proposition-Definition}

\theoremstyle{definition}
\newtheorem{defn}[thm]{Definition}
\newtheorem{exa}[thm]{Example}
\newtheorem{rem}[thm]{Remark}
\newtheorem{nota}[thm]{Notation}

\newcommand{\Span}[1]{\left<#1\right>}
\newcommand{\Hilb}[1]{#1\text{-Hilb}(\mathbb{C}^{2})}
\newcommand{\Mtheta}{\mathcal{M}_{\theta}(Q,R)}
\newcommand{\C}{\mathbb C}
\newcommand{\R}{\mathbb R}

\newcommand{\Z}{\mathbb Z}

\newcommand{\kk}{\ensuremath{\Bbbk}}

\DeclareMathOperator{\GL}{GL}
\DeclareMathOperator{\SL}{SL}

\DeclareMathOperator{\Irr}{Irr}
\DeclareMathOperator{\HILB}{Hilb}
\DeclareMathOperator{\Det}{det}
\DeclareMathOperator{\BD}{BD}
\DeclareMathOperator{\McKay}{McKayQ}
\DeclareMathOperator{\inv}{inv}
\DeclareMathOperator{\Hom}{Hom}
\DeclareMathOperator{\id}{id}

\DeclareMathOperator{\diag}{diag}
\DeclareMathOperator{\GCD}{gcd}

\DeclareMathOperator{\Sym}{Sym}

\DeclareMathOperator{\End}{End}

\title{Dihedral $G$-Hilb via representations of the McKay quiver} 
\author{\'{A}lvaro Nolla de Celis} 

\begin{document}
\maketitle

\begin{abstract}
For a given finite small binary dihedral group $G\subset\GL(2,\C)$ we provide an explicit description of the minimal resolution $Y$ of the singularity $\C^2/G$. The minimal resolution $Y$ is known to be either the moduli space of $G$-clusters $\Hilb{G}$, or the equivalent $\Mtheta$, the moduli space of $\theta$-stable quiver representations of the McKay quiver. We use both moduli approaches to give an explicit open cover of $Y$, by assigning to every distinguished $G$-graph $\Gamma$ an open set $U_\Gamma\subset\Mtheta$, and calculating the explicit equation of $U_\Gamma$ using the McKay quiver with relations $(Q,R)$.  
\end{abstract}

\tableofcontents

\section{Introduction.}

The generalisation of the McKay correspondence \cite{McK80}, \cite{Rei02} to small finite subgroups $G\subset\GL(2,\C)$ was established after Wunram \cite{Wun88} introduced the notion of special representation. The so-called ``special'' McKay correspondence relates the $G$-equivariant geometry of $\C^2$ and the minimal resolution $Y$ of the quotient $\C^2/G$, establishing a one-to-one correspondence between the irreducible components of the exceptional divisor $E\subset Y$ and the special irreducible representations. This minimal resolution $Y$ can be viewed as two equivalent moduli spaces: by a result of Ishii \cite{Ish02} it is known that $Y=\Hilb{G}$ the $G$-invariant Hilbert scheme introduced by Ito and Nakamura \cite{IN99}, and at the same time as $Y=\Mtheta$ the moduli space of $\theta$-stable representations of the McKay quiver. 

In the same spirit as \cite{Leng} in this paper we treat the problem of describing $\Hilb{G}$ by giving an explicit affine open cover. In \cite{Nak01} Nakamura introduced the notion of $G$-graphs, providing a nice and friendly framework to describe $\Hilb{G}$ for finite abelian subgroups in $\GL(n,\C)$. In this paper we consider the non-abelian analogue of a $G$-graph and provide an explicit method to interpret $\theta$-stable representations of the McKay quiver from $G$-graphs and vice versa. By using the relations on the McKay quiver, this led us to describe explicitly an open cover $\Mtheta$ (hence for $\Hilb{G}$) for binary dihedral subgroups in $\GL(2,\C)$ with the minimal number of open sets. Our method also recovers the ideals defining the $G$-clusters in $\Hilb{G}$.

An alternative description of an open cover for the minimal resolution $Y$ has been discovered independently by Wemyss \cite{WemD1}, \cite{WemD2} by using reconstruction algebras instead of the skew group ring. 

I would like to thank M. Reid for his support, the Mathematics Institute of the University of Warwick for their hospitality during my PhD, and A. Craw, Y. Ito, D. Maclagan and M. Wemyss for many useful conversations.

\section{Preliminaries.}
\subsection{Dihedral groups $\BD_{2n}(a)$ in $\GL(2,\C)$.} \label{DihedralGroups}

Let $G$ be a finite small binary dihedral subgroup in $\GL(2,\C)$. In terms of its action on the complex plane $\C^2$ we consider the representation of $G$, denoted by $\BD_{2n}(a)$, generated by 
\[ 
\text{$\alpha=
\begin{pmatrix}
\varepsilon & 0 \\
0 & \varepsilon^{a} \\
\end{pmatrix}$ and $\beta= 
\begin{pmatrix}
0 & 1 \\
-1 & 0 \\
\end{pmatrix}$}
\]
subject to relations
\[
\text{$a^2\equiv 1$ (mod $2n$), $\GCD(a+1,2n)\nmid n$}
\]
\noindent where $\varepsilon$ is a primitive $2n$-th root of unity. The group $\BD_{2n}(a)$ has order $4n$ and it contains the maximal normal index 2 cyclic subgroup $A:=\Span{\alpha}\unlhd G$, which we denote by $\frac{1}{2n}(1,a)$ (note that $\beta^2\in A$). The condition $a^2\equiv 1$ (mod $2n$) is equivalent to the relation $\alpha\beta=\beta\alpha^{a}$, and $\gcd(a+1,2n)\nmid n$ implies that the group is small (see \cite{thesis}, $\S$3 for details). 

\begin{defn}\label{defn:qk} Let $q:=\frac{2n}{(a-1,2n)}$, and $k$ such that $n=kq$. 
\end{defn} 

The group $\BD_{2n}(a)$ has $4k$ irreducible 1-dimensional representations $\rho_{j}^+$ and $\rho_j^-$ of the form
\[
\begin{array}{cc}
\rho_{j}^\pm(\alpha)=\varepsilon^j, & \rho_{j}^\pm(\beta)={\small{\left\{\begin{array}{ll}\pm i&\text{ if $n,j$ odd}\\\pm1&\text{ otherwise}\end{array}\right.}}
\end{array}
\]
where $\varepsilon$ is a $2n$-th primitive root of unity and $j$ is such that $j\equiv aj$ (mod $2n$). The values $r$ for which $r\not\equiv ar$ (mod $2n$) form in pairs the $n-k$ irreducible 2-dimensional representations $V_{r}$ of the form
\[
\begin{array}{cc}
V_{r}(\alpha)=\begin{pmatrix}\varepsilon^r&0\\0&\varepsilon^{ar}\\\end{pmatrix}, & V_{r}(\beta)=\begin{pmatrix}0&1\\(-1)^r&0\end{pmatrix}
\end{array}
\]
By definition, the natural representation is $V_{1}$. 

In what follows we take the notation as in \cite{Yos} $\S10$. Let $V(=V_1)$ a vector space with basis $\{x,y\}$ where $G$ acts naturally. Define $S=\Sym V:=\C[V^*]$ the polynomial ring in the variables $x$ and $y$. Then the action of $G$ extends to $S$ by $g\cdot f(x,y):=f(g(x),g(y))$ for $f\in S$, $g\in G$.
\begin{defn}\label{def-action} Let $G=\BD_{2n}(a),f\in S$. 
\[\begin{array}{cc}
{\small{f\in\rho_{j}^\pm :\Longleftrightarrow \alpha(f)=\varepsilon^jf, \beta(f)=\left\{\begin{array}{ll}\pm if&\text{ if $n,j$ odd}\\\pm f&\text{ otherwise}\end{array}\right.}}  \\
(f,\beta(f))\in V_k   :\Longleftrightarrow \alpha(f,\beta(f))=(\varepsilon^kf,\varepsilon^{ak}\beta(f)) 
\end{array}
\]
\end{defn}

Let $S_{\rho}:=\{f\in\C[x,y]: f\in\rho\}$ the $S^G$-module of $\rho$-invariants. Note that these are precisely the Cohen Macaulay $S^G$-modules $S_\rho=(S\otimes\rho^*)^G$ where $G$ acts on $S$ as above and $G$ acts on a representation $\rho$ by the inverse transpose. 

\subsection{$G$-Hilb and $G$-graphs.}

Let $G=\BD_{2n}(a)\subset\GL(2,\C)$ be a binary dihedral subgroup. 
\begin{defn}\label{GHilb} A $G$-{\em cluster} is a $G$-invariant zero dimensional subscheme $\mathcal{Z}\subset\C^2$ such that $\mathcal{O}_{\mathcal{Z}}\cong\C[G]$ the regular representation as $G$-modules. The $G$-{\em Hilbert scheme} $G$-Hilb($\C^2$) is the moduli space parametrising $G$-clusters. 
\end{defn}

Recall that $\C[G]=\bigoplus_{\rho\in\Irr G}(\rho)^{\dim\!\rho}$, where every irreducible representation $\rho$ appears dim$\rho$ times in the sum. Thus, as a vector space, $\mathcal{O}_{\mathcal{Z}}$ has in its basis dim$\rho$ elements in each $\rho$. To describe a distinguished basis of $\mathcal{O}_{\mathcal{Z}}$ with this property, it is convenient to use the notion of $G$-graph. 

\begin{defn}\label{defnGgraph} Let $G=\BD_{2n}(a)$. A {\em G-graph} is a subset $\Gamma\subset\C[x,y]$ satisfying the following:
\begin{enumerate}
\item It contains dim$\rho$ number of elements in each irreducible representation $\rho$. 
\item If a monomial $x^{\lambda_1}y^{\lambda_2}$ is a summand of a polynomial $P\in\Gamma$, then for every $0\leq \mu_j\leq\lambda_j$, the monomial $x^{\mu_1}y^{\mu_2}$ must be a summand of some polynomial $Q_{\mu_{1},\mu_{2}}\in\Gamma$.
\end{enumerate}
\end{defn}


For any $G$-graph $\Gamma$ there exists an open set $U_{\Gamma}\subset$ G-Hilb($\C^2$) consisting of all $G$-clusters $\mathcal{Z}$ such that $\mathcal{O}_{\mathcal{Z}}$ admits $\Gamma$ for basis as a vector space. It is proved in \cite{NdC1} (see Theorem \ref{ABCD}) that given the set of all possible $G$-graphs $\{\Gamma_i\}$, their union covers $G$-Hilb($\C^2$). 

\begin{exa}\label{exa-graph} $\Gamma=\{1,x,x^2,y,xy\}$ is a $\frac{1}{5}(1,3)$-graph. For the non-abelian binary dihedral group $D_{4}=\Span{\frac{1}{4}(1,3),\left(\begin{smallmatrix}0&1\\-1&0\end{smallmatrix}\right)}\subset\SL(2,\C)$, $\Lambda=\{1,x,y,x^2+y^2,x^2-y^2,y^3,-x^3,x^4-y^4 \}$ is a $D_{4}$-graph (note that $(x,y),(y^3,-x^3)\in V_1$). 
\end{exa}

We say that an ideal $I$ {\em represents} a $G$-graph $\Gamma$, and we write $I_\Gamma$, if $\C[x,y]/I$ admits $\Gamma$ as basis. 

In Example \ref{exa-graph}, $I_\Gamma=(x^3,x^2y,y^2)$ and similarly $I_\Lambda=(xy,x^4+y^4)$. The pictorial description of $\Gamma$ and $\Lambda$ is shown in Figure \ref{exa-graph}. Notice that for $\Lambda$ the elements $x^2+y^2\in\rho_{2}^+$ and $x^2-y^2\in\rho_{2}^-$ are described by $x^2$ and $y^2$ respectively, and the relation $x^4+y^4=0$ identifies $x^4$ and $y^4$ in $\C[x,y]/I_\Lambda$.

\begin{figure}[htbp]
\begin{center}
\Large{
\begin{pspicture}(0,0)(5,1.35)
\scalebox{0.6}{
\rput(5,-0.25){
	\scalebox{0.7}{
	\psline(0,0)(4,0)(4,0.8)(0.8,0.8)(0.8,4)(0,4)(0,0)
	\psline[fillstyle=solid,fillcolor=gray](3.2,0)(4,0)(4,0.8)(3.2,0.8)(3.2,0)
	\psline[fillstyle=solid,fillcolor=gray](0,3.2)(0,4)(0.8,4)(0.8,3.2)(0,3.2)
	\rput(0.4,0.4){$1$}
	  \rput(1.2,0.33){$x$}\rput(2,0.4){$x^2$}\rput(2.8,0.4){$x^3$}\rput(3.6,0.4){$x^4$}
	  \rput(0.4,1.2){$y$}\rput(0.4,2){$y^2$}\rput(0.4,2.8){$y^3$}\rput(0.4,3.6){$y^4$}
	  }}
\rput(0.5,-0.25){
	\scalebox{0.7}{
	\psline(0,0)(2.4,0)(2.4,0.8)(1.6,0.8)(1.6,1.6)(0,1.6)(0,0)
	\rput(0.4,0.4){$1$}
	  \rput(1.125,0.33){$x$}\rput(2,0.4){$x^2$}
	   \rput(1.125,1.2){$xy$}\rput(0.4,1.2){$y$}
	  }}
	  }
\rput(0,0.5){\normalsize $\Gamma$}\rput(2.7,1){\normalsize $\Lambda$}
\end{pspicture}	}
\caption{Representation of the $G$-graphs $\Gamma$ and $\Lambda$.}
\label{exa-graph}
\end{center}
\end{figure}
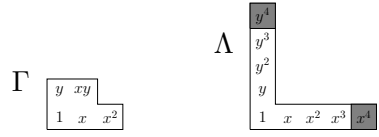

\section{$G$-graphs for $\BD_{2n}(a)$ groups.}

Let $G=\BD_{2n}(a)$. The minimal resolution $Y$ of $\C^2/G$ is obtained as follows (see \cite{IN99} $\S$1.2): First act with $A$ on $\C^2$ and consider $\Hilb{A}$ as the minimal resolution of $\C^2/A$. To complete the action of $G$ act with $G/A\cong\Span{\bar{\beta}}$ on $\Hilb{A}$. The conditions $a^2\equiv1$ (mod $2n$) and $\GCD(a+1,2n)\nmid n$ imply that the continued fraction $\frac{2n}{a}$ is symmetric with respect to the middle entry. Then the coordinates along the exceptional divisor $E=\bigcup^{2m-1}_{i=1} E_i$ are symmetric with respect to the middle curve $E_m$. The action of $G/A$ identifies the rational curves on $E$ pairwise except in $E_m$ where we have an involution. Thus the quotient $\widetilde{Y}=\Hilb{A}/(G/A)$ has two $A_1$ singularities, and the blow-up of these two points gives $\Hilb{G}$ by the uniqueness of minimal models of surfaces. 

Let us now translate this construction into graphs. Any $G/A$-orbit in $\Hilb{A}$ consists of two $A$-clusters $\mathcal{Z}$ and $\beta({\mathcal{Z}})$, with symmetric $A$-graphs $\Gamma$ and $\beta(\Gamma)$ respectively. They are represented by 
\[
\text{$I_{\Gamma}=(x^s,y^{u},x^{s-v}y^{u-r})$ and $I_{\beta(\Gamma)}=(y^{s},x^u,x^{u-r}y^{s-v})$}
\]
respectively, where $e_{i}=\frac{1}{2n}(r,s), e_{i+1}=\frac{1}{2n}(u,v)$ are two consecutive lattice points in the boundary of the Newton polygon of the lattice $L:=\Z^2+\frac{1}{2n}(1,a)\cdot\Z$. 

If we denote by $\mathcal{Y}$ the corresponding $G$-cluster, then it is clear that $\mathcal{Z}\cup\beta(\mathcal{Z})\subset\mathcal{Y}$. Thus $I_{\mathcal{Y}}\subset I_{\mathcal{Z}}\cap I_{\beta(\mathcal{Z})}$ which implies that $\widetilde{\Gamma}:=\Gamma\cup\beta(\Gamma)\subset\Gamma_{\mathcal{Y}}$. Note that the representation of $\widetilde{\Gamma}$ in the lattice of monomials is symmetric with respect to the diagonal, and the inclusion $\widetilde{\Gamma}\subset\Gamma_{\mathcal{Y}}$ is never an equality since $\Gamma$ and $\beta(\Gamma)$ always share a common subset of elements $R\subset\widetilde{\Gamma}$. The subset $\widetilde{\Gamma}$ is called $qG$-{\em graph}. 

Thus, to obtain a $G$-graph from $\widetilde{\Gamma}$ we must add $\sharp R$ elements to $\widetilde{\Gamma}$ preserving the representation spaces contained in $R$ according to Definition \ref{defnGgraph}. It is shown in \cite{NdC1} that the extension from a $qG$-graph $\widetilde{\Gamma}$ to a $G$-graph $\Gamma$ is unique. 

The following theorem resumes the classification of $G$-graphs for $\BD_{2n}(a)$ groups describing their defining ideals in each case.

\begin{thm}[\cite{NdC1}] \label{ClaGGraph} Let $G=\BD_{2n}(a)$ be a small binary dihedral group and let $\Gamma_i$ be the $A$-graph corresponding to the two consecutive lattice points $e_{i}=\frac{1}{2n}(r,s)$, $e_{i+1}=\frac{1}{2n}(u,v)$ of the Newton polygon of the lattice $L$. Denote by $\Gamma:=\Gamma(r,s;u,v)$ the $G$-graph corresponding to the $qG$-graph $\Gamma_i\cup\beta(\Gamma_i)$. 

Then we have the following possibilities:
\begin{enumerate}
\item If $u<s-v$ then $\Gamma$ is of type $A$ and it is represented by the ideal \\$I_{A}=(x^{u}y^{u},x^{s-v}y^{u-r}+(-1)^{u-r}x^{u-r}y^{s-v}, x^{r+s}+(-1)^{r}y^{r+s})$.

\item If $u-r=s-v:=m$ then $\Gamma$ is of type $B$ and 

\noindent (a) If $u<2m$ then $\Gamma$ is of type $B_1$ and it is represented by the ideal $I_{B_1}=(x^{r+s}+(-1)^{r}y^{r+s}$, $x^{m+s}y^{m-r}+(-1)^{m-r}x^{m-r}y^{m+s}$, $x^{u}y^{m}$, $x^{m}y^{u})$.

\noindent (b) If $u\geq 2m$ then $\Gamma$ is of type $B_2$ and it is represented by the ideal $I_{B_2}=(x^{2m}y^{2m}$, $x^{s+m}$, $y^{s+m}$, $x^{u}y^{m}$, $x^{m}y^{u})$.\\
	
In addition, when $u=v=q:=\frac{2n}{(a-1,2n)}$ we have four $G$-graphs of types $C^+$, $C^-$, $D^+$ and $D^-$. 

\item The $G$-graphs of types $D^\pm$ are represented by the ideals $I_{D^\pm}=(x^q\pm(-i)^qy^q, x^{s-r}y^{s-r})$.

\item For $G$-graphs of types $C^\pm$ we have two cases:

\noindent (a) If $2q<s$, and we call $m_{1}:=s-q$ and $m_{2}:=q-r$, then they are represented by the ideals $I_{C^\pm_{A}}=((x^{q}\pm(-i)^qy^{q})^2,x^sy^{m_{2}}\pm(-1)^ri^qx^{m_{2}}y^s,x^{m_{1}}y^{m_{2}}\pm(-1)^{m_{2}}x^{m_{2}}y^{m_{1}})$.

\noindent (b) If $2q=r+s$ then $I_{C^\pm_{B}}=(y^{m}(x^q\pm(-i)^qy^{q})$, $x^{m}(x^q\pm(-i)^qy^q)$,$x^{s-r}y^{s-r}$, $x^sy^{m}, x^{m}y^s)$.
\end{enumerate}
\end{thm}

\begin{rem} The list of ideals in Theorem \ref{ClaGGraph} define in $\Hilb{G}$ the intersection points of two of the exceptional curves plus the strict transform of the coordinate axis in $\C^2$.
\end{rem}

\begin{exa} Consider the $\frac{1}{12}(1,7)$-graphs given by $I_{\Gamma}=(x^7,y^2,x^{5}y)$ and $I_{\beta(\Gamma)}=(y^7,x^2,xy^5)$, with $r=1$,$s=7$,$u=2$,$v=2$. The overlap subset is $R=\{1,x,y,xy\}$ where $1\in\rho^+_0$, $xy\in\rho^-_8$ and $(x,y)\in V_1$. Then we must add the elements $x^5y-xy^5\in\rho^+_0$, $x^8-y^8\in\rho^-_8$ and $(y^7,-x^7)\in V_1$. The $\BD_{12}(5)$-graph is represented by $(x^2y^2,x^5y-xy^5,x^8-y^8)$.
\end{exa}

\begin{thm}[\cite{NdC1}] \label{ABCD} Let $G=\BD_{2n}(a)$ be small and let $P\in G$-$\HILB(\C^2)$ be defined by the ideal $I$. Then we can always choose a basis for $\C[x,y]/I$ from the list $\Gamma_{A}, \Gamma_{B}, \Gamma_{C^+}, \Gamma_{C^-}, \Gamma_{D^+}, \Gamma_{D^-}$.
Moreover, if $\Gamma_{0}, \ldots, \Gamma_{m-1}, \Gamma_{C^+}, \Gamma_{C^-}, \Gamma_{D^+}, \Gamma_{D^-}$ is the list of $G$-graphs, then an open cover of $G$-$\HILB(\C^2)$ is given by 
\[
U_{\Gamma_{0}}, \ldots, U_{\Gamma_{m-1}},  U_{\Gamma_{C^+}}, U_{\Gamma_{C^-}}, U_{\Gamma_{D^+}}, U_{\Gamma_{D^-}}.
\] 
\end{thm}

\section{$\Mtheta$ and Orbifold McKay quiver.}\label{section-quiver}

Let $G=\BD_{2n}(a)$ and let $A=\frac{1}{2n}(1,a)\unlhd G$. Denote by $\Irr G$ the set of irreducible representations of $G$. For the background material on quivers refer to \cite{ASS06}. We consider left modules (and actions), and by a path $pq$ we mean $p$ followed by $q$. Let $(Q,R)$ a quiver with relations, fix ${\bf d}=(d_i)_{i\in Q_0}$ the dimension vector of the representations of $(Q,R)$, and let $\mathbb{V}(I_R)\subset\mathbb{A}^N\cong\bigoplus_{a\in Q_1}$Mat$_{d_{t(a)}\times d_{h(a)}}$ the representation space subject to the ideal of relations $I_R$. For $\theta$ generic we define $\mathcal{M}_\theta:=\Mtheta=\mathbb{V}(I_R)/\!/\!_\theta\prod\GL(d_i)$ the moduli space of $\theta$-stable representations of $(Q,R)$ (see \cite{King}, \cite{CMT1}). Taking $Q$ to be the McKay quiver and a particular choice of generic $\theta$ (see $\S$\ref{ExplicitGHilb}) it is well known that $\mathcal{M}_\theta\cong\Hilb{G}$.

The McKay quiver of $G$ is defined by having one vertex for every $\rho\in\Irr G$ and by the number of arrows from $\rho$ to $\sigma$ to be $\dim_\C\Hom_{\C G}(\rho\otimes V,\sigma)$. Equivalently, due to Auslander it is known that the McKay quiver of $G$ is the underlying quiver of the algebra $\End_{S^G}(\bigoplus_{\rho\in\Irr G}S_\rho)$ where $S_\rho=(S\otimes\rho^*)^G$ as in \ref{def-action} (see \cite{Yos} for a proof in dimension 2).

The McKay quiver of $A$ can be drawn on a torus as follows: Let $M\cong\Z^2$ be the lattice of monomials and $M_{\inv}\cong\Z^2$ the sublattice of invariant monomials by $A$. If we take $M_{\R}=M\otimes_{\Z}\R$ we can  consider the torus $T:=M_{\R}/M_{\inv}$. The vertices are precisely $Q_{0}=M\cap T$, and the arrows between vertices are the natural multiplications by $x$ and $y$ in $M$. It is easy to see that we can always choose a fundamental domain $\mathcal{D}$ for $T$ to be the parallelogram with vertices $0$, $(k,k)$, $(2q,0)$ and $(k+2q,k)$ where the opposite sides are identified. 

\begin{figure}[h]
\begin{center}
\begin{pspicture}(0,-0.15)(17,6)
	\psset{nodesep=1pt,arcangle=15,arrowlength=2}
\scalebox{0.9}{
\rput(1.8,0.5){
	\rput(0,0){$-$}\cput(-0.5,0.5){+}
	\psdots(1.2,0)(2.4,0)(3.6,0)(4.5,0)(5.7,0)
	\rput(6.9,0){+}\rput(7.4,-0.5){$-$}
	\rput(4.1,0){$\cdots$}
	\pcline{->}(-0.35,0.42)(1.2,0)\Aput[0.1pt]{$a_{0}$}
	\pcline{->}(0.15,0)(1.2,0)\Bput[0.5pt]{$c_{0}$}
	\pcline{->}(1.2,0)(2.4,0)\Aput[0.5pt]{$r_{1,1}$}
	\pcline{->}(2.4,0)(3.6,0)\Aput[0.5pt]{$r_{1,2}$}
	\pcline{->}(4.5,0)(5.7,0)\Aput[0.5pt]{$r_{1,q-2}$}
	\pcline{->}(5.7,0)(6.7,0)\Aput[0.5pt]{$e_1$}
	\pcline{->}(5.7,0)(7.3,-0.48)\Bput[0.5pt]{$g_{1}$}
	\rput(1.2,1.2){$-$}\rput(0.7,1.7){+}
	\psdots(2.4,1.2)(3.6,1.2)(4.8,1.2)(5.7,1.2)(6.9,1.2)
	\rput(8.1,1.2){+}\rput(8.6,0.7){$-$}
	\rput(5.3,1.2){$\cdots$}
	\pcline{->}(0.85,1.62)(2.4,1.2)\Aput[0pt]{$a_{1}$}
	\pcline{->}(1.35,1.2)(2.4,1.2)\Bput[0.5pt]{$c_{1}$}
	\pcline{->}(2.4,1.2)(3.6,1.2)\Aput[0.5pt]{$r_{2,1}$}
	\pcline{->}(3.6,1.2)(4.8,1.2)\Aput[0.5pt]{$r_{2,2}$}
	\pcline{->}(5.7,1.2)(6.9,1.2)\Aput[0.5pt]{$r_{2,q-2}$}
	\pcline{->}(6.9,1.2)(7.9,1.2)\Aput[0.5pt]{$e_{2}$}
	\pcline{->}(6.9,1.2)(8.5,0.72)\Bput[0.1pt]{$g_{2}$}
	\rput(2.4,2.4){$-$}\rput(1.9,2.9){+}
	\psdots(3.6,2.4)(4.8,2.4)(6,2.4)(6.9,2.4)(8.1,2.4)
	\rput(9.3,2.4){+}\rput(9.8,1.9){$-$}
	\rput(6.5,2.4){$\cdots$}
	\pcline{->}(2.05,2.82)(3.6,2.4)\Aput[0.5pt]{$a_{2}$}
	\pcline{->}(2.55,2.4)(3.6,2.4)\Bput[0.5pt]{$c_{2}$}
	\pcline{->}(3.6,2.4)(4.8,2.4)\Aput[0.5pt]{$r_{3,1}$}
	\pcline{->}(4.8,2.4)(6,2.4)\Aput[0.5pt]{$r_{3,2}$}
	\pcline{->}(6.9,2.4)(8.1,2.4)\Aput[0.5pt]{$r_{3,q-2}$}
	\pcline{->}(8.1,2.4)(9.1,2.4)\Aput[0.5pt]{$e_{3}$}
	\pcline{->}(8.1,2.4)(9.7,1.92)\Bput[0.1pt]{$g_{3}$}

	\pcline{->}(1.2,0)(1.2,1.15)\bput[0.5pt](0.6){$d_{1}$}
	\pcline{->}(2.4,1.2)(2.4,2.35)\bput[0.5pt](0.6){$d_{2}$}

	\pcline{->}(1.2,0)(0.75,1.5)\aput[0.05pt](0.6){$b_{1}$}
	\pcline{->}(2.4,1.2)(1.95,2.7)\aput[0.05pt](0.6){$b_{2}$}

	\pcline{->}(2.4,0)(2.4,1.2)\mput*[0.5pt]{{$u_{1,1}$}}
	\pcline{->}(3.6,0)(3.6,1.2)\mput*[0.5pt]{{$u_{1,2}$}}
	\pcline{->}(5.7,0)(5.7,1.2)\mput*[0.5pt]{{$u_{1,q-2}$}}

	\pcline{->}(3.6,1.2)(3.6,2.4)\mput*[0.5pt]{$u_{2,1}$}
	\pcline{->}(4.8,1.2)(4.8,2.4)\mput*[0.5pt]{$u_{2,2}$}
	\pcline{->}(6.9,1.2)(6.9,2.4)\mput*[0.5pt]{$u_{2,q-2}$}

	\pcline{->}(6.9,0.2)(6.9,1.2)\aput[0.05pt](0.4){$f_{1}$}
	\pcline{->}(8.1,1.4)(8.1,2.4)\aput[0.05pt](0.4){$f_{2}$}

	\pcline{->}(7.37,-0.4)(6.9,1.2)\bput[0.05pt](0.4){$h_{1}$}
	\pcline{->}(8.57,0.8)(8.1,2.4)\bput[0.05pt](0.4){$h_{2}$}
	\rput(-0.5,1.1){$\boldsymbol{S_{\rho_{0}}}$}
	\rput(7.8,-0.7){$\boldsymbol{S_{\rho_{q}}}$}
	\rput(0.7,2.1){$\boldsymbol{S_{\rho_{a+1}}}$}
	\rput(9,0.4){$\boldsymbol{S_{\rho_{a+1+q}}}$}
	\rput(1.9,3.3){$\boldsymbol{S_{\rho_{2(a+1)}}}$}
	\rput(10,1.6){$\boldsymbol{S_{\rho_{2(a+1)+q}}}$}
	}
	\psline[linestyle=dotted](5.8,3.5)(6.3,4)
	\psline[linestyle=dotted](7.8,3.5)(8.3,4)
	\psline[linestyle=dotted](9.8,3.5)(10.3,4)
\rput(3.5,3.7){
	\rput(1.2,1){$-$}\rput(0.7,1.5){+}
	\psdots(2.4,1)(3.6,1)(4.8,1)(5.7,1)(6.9,1)
	\rput(8.1,1){+}\rput(8.6,0.5){$-$}
	\rput(5.3,1){$\cdots$}
	\pcline{->}(0.85,1.42)(2.4,1)\rput(1.6,1.4){$a_{k-1}$}
	\pcline{->}(1.35,1)(2.4,1)\Bput[0.5pt]{$c_{k-1}$}
	\pcline{->}(2.4,1)(3.6,1)\Aput[0.5pt]{$r_{0,1}$}
	\pcline{->}(3.6,1)(4.8,1)\Aput[0.5pt]{$r_{0,2}$}
	\pcline{->}(5.7,1)(6.9,1)\Aput[0.5pt]{$r_{0,q-2}$}
	\pcline{->}(6.9,1)(7.9,1)\Aput[0.5pt]{$e_{0}$}
	\pcline{->}(6.9,1)(8.5,0.52)\Bput[0.1pt]{$g_{0}$}
	\rput(2.4,2.2){$-$}\cput(1.9,2.7){+}
	\psdots(3.6,2.2)(4.8,2.2)(6,2.2)(6.9,2.2)(8.1,2.2)
	\rput(9.3,2.2){+}\rput(9.8,1.7){$-$}
	\rput(6.5,2.2){$\cdots$}
	\pcline{->}(2.05,2.62)(3.6,2.2)\Aput[0.5pt]{$a_{0}$}
	\pcline{->}(2.55,2.2)(3.6,2.2)\Bput[0.5pt]{$c_{0}$}
	\pcline{->}(3.6,2.2)(4.8,2.2)\Aput[0.5pt]{$r_{1,1}$}
	\pcline{->}(4.8,2.2)(6,2.2)\Aput[0.5pt]{$r_{1,2}$}
	\pcline{->}(6.9,2.2)(8.1,2.2)\Aput[0.5pt]{$r_{1,q-2}$}
	\pcline{->}(8.1,2.2)(9.1,2.2)\Aput[0.5pt]{$e_{1}$}
	\pcline{->}(8.1,2.2)(9.7,1.72)\Bput[0.1pt]{$g_{1}$}
	\pcline{->}(2.4,1)(1.95,2.5)\aput[0.05pt](0.6){$b_{0}$}
	\pcline{->}(2.4,1)(2.4,2.15)\bput[0.5pt](0.6){$d_{0}$}
	\pcline{->}(3.6,1)(3.6,2.2)\mput*[0.5pt]{$u_{0,1}$}
	\pcline{->}(4.8,1)(4.8,2.2)\mput*[0.5pt]{$u_{0,2}$}
	\pcline{->}(6.9,1)(6.9,2.2)\mput*[0.5pt]{$u_{0,q-2}$}
	\pcline{->}(8.1,1)(8.1,2.2)\aput[0.05pt](0.4){$f_{0}$}
	\pcline{->}(8.57,0.6)(8.1,2.2)\bput[0.05pt](0.4){$h_{0}$}
	\rput(0.3,1.9){$\boldsymbol{S_{\rho_{(k-1)(a+1)}}}$}
	\rput(9.25,0.2){$\boldsymbol{S_{\rho_{(k-1)(a+1)+q}}}$}
	\rput(1.3,2.6){$\boldsymbol{S_{\rho_{0}}}$}
	\rput(9.75,1.4){$\boldsymbol{S_{\rho_{q}}}$}
	}}
\end{pspicture}
\caption{McKay quiver for $\BD_{2n}(a)$ groups}
\label{McKayQ}
\end{center}
\end{figure}
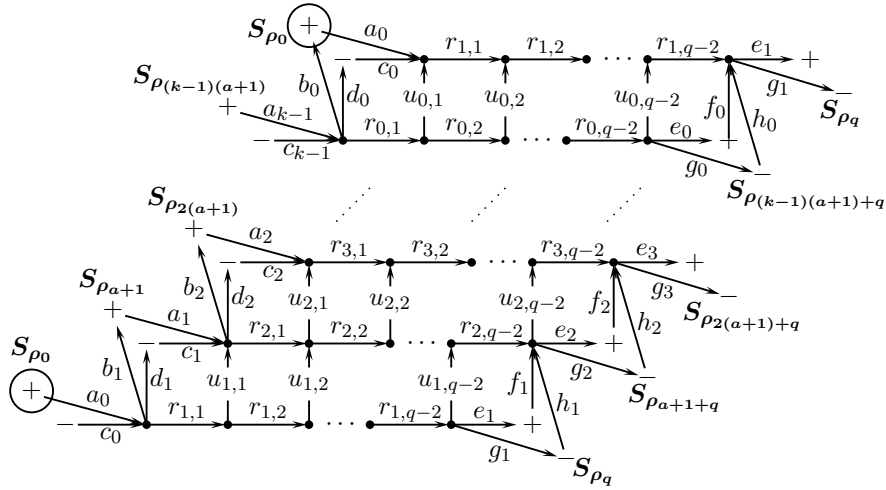

\begin{prop}\label{prop-McKayQ} (i) The McKay quiver $Q$ of $\BD_{2n}(a)$ is the $\Z/2$-orbifold quotient of the McKay quiver for the Abelian subgroup $A$ (see Figure \ref{McKayQ}). \\
(ii) The relations $R$ on the $Q$ which gives the identification between $\Hilb{G}$ and $\mathcal{M}_\theta$ are 
\begin{align*}
&a_{i}b_{i+1}=0, f_{i}e_{i+1}=0, c_{i}d_{i+1}=0, h_{i}g_{i+1}=0, \\
&b_{i}a_{i}+d_{i}c_{i}=r_{i,1}u_{i,1}, e_{i}f_{i}+g_{i}h_{i}=u_{i,q-2}r_{i+1,q-2}, \\
&u_{i,j}r_{i+1,j}=r_{i,j+1}u_{i,j+1},
\end{align*} 
considering subindices modulo $k$.
\end{prop}

\begin{nota} The source and target for $r_{i,j}$ and $u_{i,j}$ are 
\begin{align*}
r_{i,j}:&S_{V_{\overline{(i-1)}(a+1)+j}}\to S_{V_{\overline{(i-1)}(a+1)+j+1}} \\ 
u_{i,j}:&S_{V_{\overline{(i-1)}(a+1)}+j+1}\to S_{V_{i(a+1)+j}}
\end{align*}
with $i\in[0,k-1]$, $j\in[1,q-2]$, where $\overline{i}$ denotes $i$ mod $k$.
\end{nota}

\begin{rem} In the case $q=2$ the relations are $a_{i}b_{i+1}=0$, $f_{i}e_{i+1}=0$, $c_{i}d_{i+1}=0$, $h_{i}g_{i+1}=0$ and $b_{i}a_{i}+d_{i}c_{i}=e_{i}f_{i}+g_{i}h_{i}$.
\end{rem}

\begin{proof} (i) Let $\Irr A=\{\rho_{0},\ldots,\rho_{2n-1}\}$. The group $G$ acts on $A$ by conjugation, which induces an action of $G/A\cong\Z/2$ on $A$ by $\beta\cdot h := \beta h\beta^{-1}$, for any $h\in A$. Therefore $G/A$ acts on $\Irr A$ by $\beta\cdot\rho_{k}:=\rho_{ak}$, for $\rho_{k}\in\Irr A$. The free orbits are $\{\rho_{i},\rho_{ai}\}$ with $ai\not\equiv i$ mod $2n$, producing the 2-dimensional representations $V_{i}$ in the McKay quiver of $G$, which we denote by $\McKay(G)$. Every fixed point $\rho_{j}$ with $aj\equiv j$ (mod $2n$) splits into the two 1-dimensional representations $\rho_{j}^+$ and $\rho_{j}^-$ in $\McKay(G)$. 

Note that $\McKay(G)$ is now drawn on a cylinder where only the top and bottom sides are identified. The arrows of $\McKay(A)$ going in and out fixed representations split into two different arrows, while for the rest we have a 1-to-1 correspondence between arrows in $\McKay(A)$ and $\McKay(G)$. 

(ii) Let $\kk Q$ be the path algebra, $S=\C[x,y]$ and $V^*$ the natural representation. Tensoring with $\Det_{V^*}:=\bigwedge^2V^*$ induces a permutation $\tau$ on $Q_0$ by $e_{i}=\tau({e_{j}}) \Longleftrightarrow  \rho_{i}=\rho_{j}\otimes\Det_{V^*}$. Now consider an arrow $a:e_{i}\to e_{j}$ as an element $\psi_{a}\in\Hom_{\C G}(\rho_{i},\rho_{j}\otimes V^*)$. Then for any path $p=a_{1}a_{2}$
of length 2 we can consider the $G$-module homomorphism
$\rho_{t(p)}\stackrel{\psi_{p}}{\longrightarrow} \rho_{h(p)}\otimes V^{*\otimes r} 
		\stackrel{\id_{\rho_{h(p)}}\otimes\gamma}{\longrightarrow}  \rho_{h(p)}\otimes\Det_{V^*}$
where $\psi_{p}$ is the composition of the maps $\psi_{a_1}$ and $\psi_{a_{2}}\otimes\id_V^*$, and $\gamma:V^{*\otimes2}\to\bigwedge^2V^*$ sends $v_{1}\otimes v_{2}\mapsto v_{1}\wedge v_{2}$. By Schur's Lemma the composition of the maps above is zero if $\tau(h(p))\neq t(p)$, 
a scalar $c_{p}$ otherwise. It is known by \cite{BSW} that for a finite small $G\subset\GL(2,\C)$ (and more generally for any small finite subgroup $G\subset\GL(r,\C)$) the skew group algebra $S\!\ast\!G$ is Morita equivalent to the algebra $\kk Q/\Span{\partial_{p}\Phi : |p|=0}$, where $\Phi := \sum_{|p|=2}(c_{p}\dim h(p))p$ and $\partial_{p}$ are derivations with respect to paths of length $0$, i.e.\ vertices $e_i$. Since the $\theta$-stable $S\ast G$-modules ($\theta$ as in $\S$\ref{ExplicitGHilb}) are precisely the $G$-clusters, which gives $\mathcal{M}_\theta\cong G$-Hilb($\C^2$). 

For $G=\BD_{2n}(a)$, $\Det_{V^*}=\rho^+_{a+1}$ so $\tau$ translates $\McKay(G)$ one step diagonally up (see Figure \ref{McKayQ}). Only paths of length 2 joining two vertices identified by $\tau$ appear in $\Phi$, giving the relations $R$ by derivations with respect to the vertices of $Q$. 
\end{proof}

\begin{exa}\label{Exa:ReprBD12} Consider the group $\BD_{30}(19)$ generated by $\alpha=\diag(\varepsilon,\varepsilon^{19})$ with $\varepsilon$ a primitive 30-$th$ root of unity, and $\beta=\left(\begin{smallmatrix}0&1\\-1&0\end{smallmatrix}\right)$. We have $q=5$ and $k=3$. The continued fraction $\frac{30}{30-19}=\frac{30}{11}=[3,4,3]$ describes the lattice $M_{\inv}$. The two consecutive invariant monomials $x^3y^3$ and $x^{11}y$ define a fundamental domain of 
the lattice $T$, which can be translated into the parallelogram filled with numbers shown 
in Figure \ref{McKayBD30(19)} (a). The diagram represents the lattice $M$ where the bottom left corner represents the monomial $1$ and the numbers denote the representation to which they belong to, e.g. the number \textbf{0} corresponds to monomials in $M_{\inv}$. Opposite sides of the parallelogram are identified. The McKay quiver is completed by adding at every vertex the two arrows corresponding to the multiplication by $x$ and $y$ to the corresponding adjacent vertices. 

Now acting by $\beta$ we see that representations $\rho_{0}, \rho_{20}, \rho_{10}$ and $\rho_{5}, \rho_{25}, \rho_{15}$ are fixed, while the rest (in pairs) are contained in a free orbit. 
The McKay quiver for $\BD_{30}(19)$ is shown in Figure \ref{McKayBD30(19)} (b). Notice that top and bottom rows are identified.  

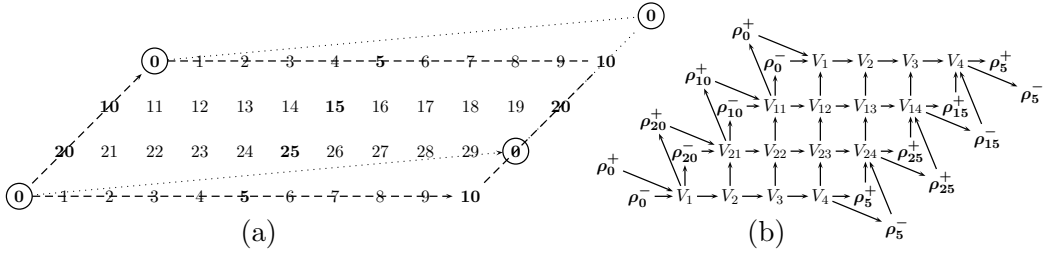
\begin{figure}[htbp]
\begin{center}
\begin{pspicture}(0,-0.5)(18,2.5)
\rput(0,0){
\scalebox{0.6}{
	\psset{nodesep=5pt}
	\cput(0,0){\rnode{P1}{$\boldsymbol{0}$}}
	\rput(1,0){1}\rput(2,0){2}\rput(3,0){3}\rput(4,0){4}\rput(5,0){$\boldsymbol{5}$}
	\rput(6,0){6}\rput(7,0){7}\rput(8,0){8}\rput(9,0){9}
	\rput(10,0){\rnode{P5}{$\boldsymbol{10}$}}
	
	\rput(1,1){$\boldsymbol{20}$}\rput(2,1){21}\rput(3,1){22}\rput(4,1){23}\rput(5,1){24}
	\rput(6,1){$\boldsymbol{25}$}\rput(7,1){26}\rput(8,1){27}\rput(9,1){28}\rput(10,1){29}	
	\cput(11,1){\rnode{P3}{$\boldsymbol{0}$}}
	
	\rput(2,2){$\boldsymbol{10}$}\rput(3,2){11}\rput(4,2){12}\rput(5,2){13}\rput(6,2){14}
	\rput(7,2){$\boldsymbol{15}$}\rput(8,2){16}\rput(9,2){17}\rput(10,2){18}\rput(11,2){19}
	\rput(12,2){$\boldsymbol{20}$}
	
	\cput(3,3){\rnode{P2}{$\boldsymbol{0}$}}\rput(4,3){1}\rput(5,3){2}\rput(6,3){3}\rput(7,3){4}
	\rput(8,3){$\boldsymbol{5}$}\rput(9,3){6}\rput(10,3){7}\rput(11,3){8}\rput(12,3){9}
	\rput(13,3){\rnode{P6}{$\boldsymbol{10}$}}
	
	\cput(14,4){\rnode{P4}{$\boldsymbol{0}$}}
	\ncline[linestyle=dotted]{->}{P1}{P3}\ncline[linestyle=dashed]{->}{P1}{P2}
	\ncline[linestyle=dotted]{P2}{P4}\ncline[linestyle=dotted]{P3}{P4}
	\ncline[linestyle=dashed]{->}{P1}{P5}\ncline[linestyle=dashed]{P2}{P6}
	\ncline[linestyle=dashed]{P5}{P6}
	}}
\rput(8.25,0){
\scalebox{0.6}{
	\psset{nodesep=2pt}
	\rput(-0.7,0.7){\rnode{P0+}{$\boldsymbol{\rho^+_{0}}$}}
	\rput(0,0){\rnode{P0-}{$\boldsymbol{\rho^-_{0}}$}}
	\rput(1,0){\rnode{V1}{$V_{1}$}}\rput(2,0){\rnode{V2}{$V_{2}$}}
	\rput(3,0){\rnode{V3}{$V_{3}$}}\rput(4,0){\rnode{V4}{$V_{4}$}}
	\rput(5,0){\rnode{P5+}{$\boldsymbol{\rho^+_{5}}$}}
	\rput(5.7,-0.7){\rnode{P5-}{$\boldsymbol{\rho^-_{5}}$}}
	
	\ncline{->}{P0+}{V1}\ncline{->}{P0-}{V1}
	\ncline{->}{V1}{V2}\ncline{->}{V2}{V3}\ncline{->}{V3}{V4}
	\ncline{->}{V4}{P5+}\ncline{->}{V4}{P5-}
	
	\rput(0.3,1.7){\rnode{P20+}{$\boldsymbol{\rho^+_{20}}$}}
	\rput(1,1){\rnode{P20-}{$\boldsymbol{\rho^-_{20}}$}}
	\rput(2,1){\rnode{V21}{$V_{21}$}}\rput(3,1){\rnode{V22}{$V_{22}$}}
	\rput(4,1){\rnode{V23}{$V_{23}$}}\rput(5,1){\rnode{V24}{$V_{24}$}}
	\rput(6,1){\rnode{P25+}{$\boldsymbol{\rho^+_{25}}$}}
	\rput(6.7,0.3){\rnode{P25-}{$\boldsymbol{\rho^+_{25}}$}}

	\ncline{->}{P20+}{V21}\ncline{->}{P20-}{V21}
	\ncline{->}{V21}{V22}\ncline{->}{V22}{V23}\ncline{->}{V23}{V24}
	\ncline{->}{V24}{P25+}\ncline{->}{V24}{P25-}
	
	\rput(1.3,2.7){\rnode{P10+}{$\boldsymbol{\rho^+_{10}}$}}	\rput(2,2){\rnode{P10-}{$\boldsymbol{\rho^-_{10}}$}}
	\rput(3,2){\rnode{V11}{$V_{11}$}}\rput(4,2){\rnode{V12}{$V_{12}$}}
	\rput(5,2){\rnode{V13}{$V_{13}$}}\rput(6,2){\rnode{V14}{$V_{14}$}}
	\rput(7,2){\rnode{P15+}{$\boldsymbol{\rho^+_{15}}$}}
	\rput(7.7,1.3){\rnode{P15-}{$\boldsymbol{\rho^-_{15}}$}}

	\ncline{->}{P10+}{V11}\ncline{->}{P10-}{V11}
	\ncline{->}{V11}{V12}\ncline{->}{V12}{V13}\ncline{->}{V13}{V14}
	\ncline{->}{V14}{P15+}\ncline{->}{V14}{P15-}
	
	\rput(2.3,3.7){\rnode{P00+}{$\boldsymbol{\rho^+_{0}}$}}	\rput(3,3){\rnode{P00-}{$\boldsymbol{\rho^-_{0}}$}}
	\rput(4,3){\rnode{VV1}{$V_{1}$}}\rput(5,3){\rnode{VV2}{$V_{2}$}}
	\rput(6,3){\rnode{VV3}{$V_{3}$}}\rput(7,3){\rnode{VV4}{$V_{4}$}}
	\rput(8,3){\rnode{P55+}{$\boldsymbol{\rho^+_{5}}$}}
	\rput(8.7,2.3){\rnode{P55-}{$\boldsymbol{\rho^-_{5}}$}}
	
	\ncline{->}{P00+}{VV1}\ncline{->}{P00-}{VV1}
	\ncline{->}{VV1}{VV2}\ncline{->}{VV2}{VV3}\ncline{->}{VV3}{VV4}
	\ncline{->}{VV4}{P55+}\ncline{->}{VV4}{P55-}

	\ncline{->}{V1}{P20+}\ncline{->}{V1}{P20-}
	\ncline{->}{V2}{V21}\ncline{->}{V3}{V22}
	\ncline{->}{V4}{V23}\ncline{->}{P5+}{V24}
	\ncline{->}{P5-}{V24}

	\ncline{->}{V21}{P10+}\ncline{->}{V21}{P10-}
	\ncline{->}{V22}{V11}\ncline{->}{V23}{V12}
	\ncline{->}{V24}{V13}\ncline{->}{P25+}{V14}
	\ncline{->}{P25-}{V14}

	\ncline{->}{V11}{P00+}\ncline{->}{V11}{P00-}
	\ncline{->}{V12}{VV1}\ncline{->}{V13}{VV2}
	\ncline{->}{V14}{VV3}\ncline{->}{P15+}{VV4}
	\ncline{->}{P15-}{VV4}
	}}
\rput(3.25,-0.5){(a)}\rput(10,-0.5){(b)}
\end{pspicture}
\caption{McKay quivers for (a) the abelian group $\frac{1}{30}(1,19)$ and for (b) the group $\BD_{30}(19)$.}
\label{McKayBD30(19)}
\end{center}
\end{figure}
\end{exa}

\section{Explicit calculation of $\Hilb{G}$} \label{ExplicitGHilb}
 
Let $G=\BD_{2n}(a)$ and $(Q,R)$ be the McKay quiver as in \ref{section-quiver}. Denote the arrows by {\bf a} $=(a,A), {\bf a}=\left(\begin{smallmatrix}a\\A\end{smallmatrix}\right)$, or {\bf a} $=\left(\begin{smallmatrix}a&A\\a'&A'\end{smallmatrix}\right)$ depending on the dimensions at the source and target of {\bf a}. We take representations of $Q$ with dimension vector $\textbf{d}=($dim$\rho_{i})_{i\in Q_{0}}$ and the generic stability condition $\theta=(1-\sum_{\rho_{i}\in\Irr G}$dim$\rho_{i},1\ldots,1)$, which ensures that $\Hilb{G}\cong\mathcal{M}_\theta$. 

\begin{cla}\label{stab-cond} A representation of $Q$ is stable if and only if there exist $\dim\rho_i$ linearly independent maps from the distinguished source, chosen to be $\rho_{0}\in Q_{0}$, to every other vertex $\rho_i$ in $Q$. 
\end{cla}

Indeed, a representation $W$ is not $\theta$-stable if and only if there exists a proper subrepresentation $W'\subset W$ with $\theta(W')<0$. Since the only nonzero entry of $\theta$ corresponds to $\rho_0^+$ and $\textbf{d}(W')$ has to be strictly smaller than $\textbf{d}(W)$, this is equivalent to say that there are strictly less linearly independent paths from $\rho_0^+$ to $\rho_i$ than $\dim\rho_i$, for $i\neq0$. \\

By making a correspondence between elements of a $G$-graph and paths in $Q$, an open cover of $\mathcal{M}_\theta$ is given by the ones corresponding to the $G$-graphs. The $G$-graphs predetermine the choices of linear independent paths, thus giving a covering of $\mathcal{M}_\theta$ with the minimal number of open sets. 

Let us illustrate this fact by a simple example.

\begin{exa}\label{CyclicExa} Consider the group \[G=\frac{1}{5}(1,2)=\Span{\alpha=\left(\begin{smallmatrix}\varepsilon&0\\0&\varepsilon^2\end{smallmatrix}\right)|\varepsilon^5=1, \text{$\varepsilon$~primitive}}.\] There are five 1-dimensional irreducible representations $\rho_{i}(\alpha)=\varepsilon^{i}$ for $i=0,\ldots,4$. The McKay quiver with relations $(Q,R)$ and the corresponding quiver between the $\C[x,y]^G$-modules $S_{\rho_{i}}$ are shown in Figure \ref{MQuivZ5(1,2)}. 

\begin{figure}[htbp]
\begin{center}
\begin{pspicture}(0,-0.25)(8,2.5)
\scalebox{0.9}{
	\psset{nodesep=3pt,arcangle=-20}
	\rput(0,0){
	\rput(1.5,2.5){\rnode{0}{$\rho_{0}$}}
	\rput(0,1.5){\rnode{1}{$\rho_{1}$}}
	\rput(0.5,0){\rnode{2}{$\rho_{2}$}}
	\rput(2.5,0){\rnode{3}{$\rho_{3}$}}
	\rput(3,1.5){\rnode{4}{$\rho_{4}$}}
	\ncarc{->}{0}{2}\rput(0.9,1.25){A}\ncline{->}{0}{1}\Bput*[0.15]{a}
	\ncarc{->}{1}{3}\rput(1.2,0.7){B}\ncline{->}{1}{2}\Bput*[0.15]{b}
	\ncarc{->}{2}{4}\rput(1.8,0.7){C}\ncline{->}{2}{3}\Bput*[0.15]{c}
	\ncarc{->}{3}{0}\rput(2.1,1.25){D}\ncline{->}{3}{4}\Bput*[0.15]{d}
	\ncarc{->}{4}{1}\rput(1.5,1.6){E}\ncline{->}{4}{0}\Bput*[0.15]{e}
	\rput(5,2){aB = Ac}\rput(5,1.5){bC = Bd}
	\rput(5,1){cD = Ce}\rput(5,0.5){dE = Da}
	\rput(5,0){eA = Eb}
	}
	\rput(7,0){
	\rput(1.5,2.5){\rnode{0}{$S_{\rho_{0}}$}}
	\rput(0,1.5){\rnode{1}{$S_{\rho_{1}}$}}
	\rput(0.5,0){\rnode{2}{$S_{\rho_{2}}$}}
	\rput(2.5,0){\rnode{3}{$S_{\rho_{3}}$}}
	\rput(3,1.5){\rnode{4}{$S_{\rho_{4}}$}}
	\ncarc{->}{0}{2}\rput(0.9,1.25){$y$}\ncline{->}{0}{1}\Bput*[0.15]{$x$}
	\ncarc{->}{1}{3}\rput(1.2,0.7){$y$}\ncline{->}{1}{2}\Bput*[0.15]{$x$}
	\ncarc{->}{2}{4}\rput(1.8,0.7){$y$}\ncline{->}{2}{3}\Bput*[0.15]{$x$}
	\ncarc{->}{3}{0}\rput(2.1,1.25){$y$}\ncline{->}{3}{4}\Bput*[0.15]{$x$}
	\ncarc{->}{4}{1}\rput(1.5,1.6){$y$}\ncline{->}{4}{0}\Bput*[0.15]{$x$}
	}}
\end{pspicture}
\caption{The McKay quiver for $\frac{1}{5}(1,2)$ with relations and the quiver between the modules $S_{\rho}$.}
\label{MQuivZ5(1,2)}
\end{center}
\end{figure}
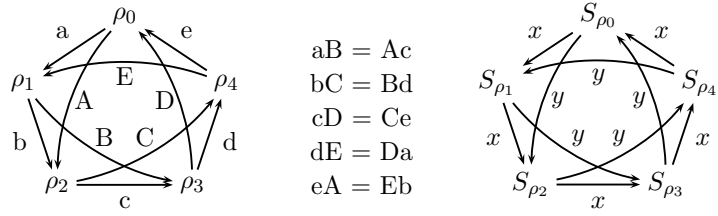

The quotient $\C^2/G$ is toric and by the continued fraction $\frac{5}{2}=[3,2]$ we know that the minimal resolution $\Hilb{\frac{1}{5}(1,2)}\cong\Mtheta$ is the union of 3 complex planes: \[U_0\cup U_1\cup U_2\cong\C^2_{\Span{A=\frac{y}{x^2},e=x^5}}\cup\C^2_{\Span{d=\frac{x^2}{y},E=\frac{y^3}{x}}}\cup\C^2_{\Span{a=\frac{x}{y^3},D=y^5}}.\] 

The $G$-graphs for every open set and their corresponding open sets in $\Mtheta$ are shown in Figure \ref{Z5(1,2)}. For every open set $U_i, i=0,1,2$, we also give their representation space in $\Mtheta$, and the ideals of every point in $\Hilb{\frac{1}{5}(1,2)}$. Note that multiplication by $x$ and $y$ among elements in the $G$-graphs determine the linearly independent (nonzero in this case) maps from $\rho_0$ of Claim \ref{stab-cond}.

\begin{figure}[htbp]
\begin{center}
\begin{pspicture}(0,1)(10,11)
	\psset{nodesep=1pt,arcangle=-20}
	\rput(2.5,10.5){$U_0\cong\C^2_{A,e}$}
	\rput(6,10.5){$U_1\cong\C^2_{d,E}$}
	\rput(9.5,10.5){$U_2\cong\C^2_{a,D}$}
	
	\rput(-0.5,9.75){\small Ideals in}
	\rput(-0.5,9.25){\small $U_i\subset\Hilb{\frac{1}{5}(1,2)}$}
	\rput(2.5,9.5){\scalebox{0.8}{$(x^5-e,y-Ax^2)$}}
	\rput(6,9.5){\scalebox{0.8}{$(x^2-dy,y^3-Ex,xy^2-bE)$}}
	\rput(9.5,9.5){\scalebox{0.8}{$(x-ay^3,y^5-D)$}}
	
	\rput(-0.5,7.5){\small $G$-graphs}
	\rput(-0.5,7){\small $\Gamma_i$}
	\rput(1,6.9){
	\scalebox{0.8}{
	\psline[linewidth=1pt](0.15,0.15)(3.4,0.15)(3.4,0.8)(0.15,0.8)(0.15,0.15)
	\rput(0.45,0.45){\rnode{1}{\red $1^{\empty}$}}
	\rput(1.05,0.44){\rnode{x}{\red $x$}}
	\rput(1.7,0.5){\rnode{x2}{\red $x^2$}}
	\rput(2.4,0.5){\rnode{x3}{\red $x^3$}}
	\rput(0.45,1){\rnode{y}{$y$}}
	\rput(3.1,0.5){\rnode{x4}{\red $x^4$}}
	\rput(3.8,0.5){\rnode{x5}{$x^5$}}
	\ncLine[linewidth=0.5pt]{->}{1}{x}
	\psline[linewidth=0.5pt]{->}(1.175,0.44)(1.475,0.435)
	\ncLine[offset=-2pt,linewidth=0.5pt]{->}{x2}{x3}
	\ncLine[offset=-2pt,linewidth=0.5pt]{->}{x3}{x4}
	}}
	\rput(5.25,6.4){
	\scalebox{0.8}{
	\psline[linewidth=1pt](0.15,0.15)(1.4,0.15)(1.4,1.4)(0.75,1.4)(0.75,2.1)(0.15,2.1)(0.15,0.15)
	\rput(0.45,0.45){\rnode{1}{\red 1}}
	\rput(1.05,0.43){\rnode{x}{\red $x$}}
	\rput(1.7,0.475){\rnode{x2}{$x^2$}}
	\rput(0.45,1){\rnode{y}{\red $y$}}
	\rput(0.5,1.7){\rnode{y2}{\red $y^2$}}
	\rput(0.5,2.4){\rnode{y3}{$y^3$}}
	\rput(1.05,1){\rnode{xy}{\red $xy$}}
	\rput(1.1,1.7){\rnode{xy2}{$xy^2$}}
	\ncline[linewidth=0.5pt]{->}{1}{x}
	\ncline[linewidth=0.5pt]{->}{1}{y}
	\psline[linewidth=0.5pt]{->}(0.45,1.175)(0.45,1.475)
	\ncline[linewidth=0.5pt]{->}{y}{xy}
	\ncline[linewidth=0.5pt,linestyle=dotted,dotsep=0.5pt]{->}{x}{xy}
	}}
	\rput(9,5.9){
	\scalebox{0.8}{
	\psline[linewidth=1pt](0.15,0.15)(0.8,0.15)(0.8,3.4)(0.15,3.4)(0.15,0.15)
	\rput(0.45,0.45){\rnode{1}{\red 1}}
	\rput(0.45,1){\rnode{y}{\red $y$}}
	\rput(0.5,1.7){\rnode{y2}{\red $y^2$}}
	\rput(0.5,2.4){\rnode{y3}{\red $y^3$}}
	\rput(0.5,3.1){\rnode{y4}{\red $y^4$}}
	\rput(0.5,3.75){\rnode{y5}{$y^5$}}
	\rput(1.05,0.43){\rnode{x}{$x$}}
	\ncline[linewidth=0.5pt]{->}{1}{y}
	\psline[linewidth=0.5pt]{->}(0.45,1.175)(0.45,1.475)
	\ncline[offset=2pt,linewidth=0.5pt]{->}{y2}{y3}
	\ncline[offset=2pt,linewidth=0.5pt]{->}{y3}{y4}
	}}
	
	\rput(-0.5,5){\small Open set}
	\rput(-0.5,4.5){\small $U_i\subset\Mtheta$}
	\rput(1.5,4){
	\scalebox{0.6}{
	\rput(1.5,2.5){\rnode{0}{\red $1$}}
	\rput(0,1.5){\rnode{1}{\red $x$}}
	\rput(0.5,0){\rnode{2}{\red $x^2$}}
	\rput(2.5,0){\rnode{3}{\red $x^3$}}
	\rput(3,1.5){\rnode{4}{\red $x^4$}}
	\ncline[linewidth=2pt]{->}{0}{1}\Bput*[0.15]{$1$}
	\ncline[linewidth=2pt]{->}{1}{2}\Bput*[0.15]{$1$}
	\ncline[linewidth=2pt]{->}{2}{3}\Bput*[0.15]{$1$}
	\ncline[linewidth=2pt]{->}{3}{4}\Bput*[0.15]{$1$}
	\rput(1.5,-1){\Large $a, b, c, d\neq0$}
	}}
	\rput(5,4){
	\scalebox{0.6}{
	\rput(1.5,2.5){\rnode{0}{\red $1$}}
	\rput(0,1.5){\rnode{1}{\red $x$}}
	\rput(0.5,0){\rnode{2}{\red $y$}}
	\rput(2.5,0){\rnode{3}{\red $xy$}}
	\rput(3,1.5){\rnode{4}{\red $y^2$}}
	\ncarc[linewidth=2pt]{->}{0}{2}\rput(0.9,1.25){$1$}
	\ncline[linewidth=2pt]{->}{0}{1}\Bput*[0.15]{$1$}
	\ncarc[linewidth=2pt]{->}{2}{4}\rput(1.8,0.7){$1$}
	\ncline[linewidth=2pt]{->}{2}{3}\Bput*[0.15]{$1$}
	\ncarc[linewidth=1pt,linestyle=dotted,dotsep=0.5pt]{->}{1}{3}\Aput*[0]{$1$}\ncarc[linewidth=1pt,linestyle=dotted,dotsep=0.5pt]{->}{1}{3}
	\rput(1.5,-1){\Large $a, A, c, C\neq0$}
	}}
	\rput(8.5,4){
	 \scalebox{0.6}{
	\rput(1.5,2.5){\rnode{0}{\red $1$}}
	\rput(0,1.5){\rnode{1}{\red $y^3$}}
	\rput(0.5,0){\rnode{2}{\red $y$}}
	\rput(2.5,0){\rnode{3}{\red $y^4$}}
	\rput(3,1.5){\rnode{4}{\red $y^2$}}
	\ncarc[linewidth=2pt]{->}{0}{2}\rput(0.9,1.25){$1$}
	\ncarc[linewidth=2pt]{->}{1}{3}\rput(1.2,0.7){$1$}
	\ncarc[linewidth=2pt]{->}{2}{4}\rput(1.8,0.7){$1$}
	\ncarc[linewidth=2pt]{->}{4}{1}\rput(1.5,1.6){$1$}
	\rput(1.5,-1){\Large $A, B, C, E\neq0$}
	}}
	
	\rput(-0.5,2){\small Representation}
	\rput(-0.5,1.5){\small space}
	\rput(1.5,1){
	 \scalebox{0.6}{
	\rput(1.5,2.5){\rnode{0}{$\C$}}
	\rput(0,1.5){\rnode{1}{$\C$}}
	\rput(0.5,0){\rnode{2}{$\C$}}
	\rput(2.5,0){\rnode{3}{$\C$}}
	\rput(3,1.5){\rnode{4}{$\C$}}
	\ncarc{->}{0}{2}\rput(0.9,1.25){$A$}\ncline{->}{0}{1}\Bput*[0.15]{$1$}
	\ncarc{->}{1}{3}\rput(1.2,0.7){$A$}\ncline{->}{1}{2}\Bput*[0.15]{$1$}
	\ncarc{->}{2}{4}\rput(1.8,0.7){$A$}\ncline{->}{2}{3}\Bput*[0.15]{$1$}
	\ncarc{->}{3}{0}\rput(2.05,1.25){$Ae$}\ncline{->}{3}{4}\Bput*[0.15]{$1$}
	\ncarc{->}{4}{1}\rput(1.5,1.6){$Ae$}\ncline{->}{4}{0}\Bput*[0.15]{$e$}
	}}
	\rput(5,1){
	 \scalebox{0.6}{
	\rput(1.5,2.5){\rnode{0}{$\C$}}
	\rput(0,1.5){\rnode{1}{$\C$}}
	\rput(0.5,0){\rnode{2}{$\C$}}
	\rput(2.5,0){\rnode{3}{$\C$}}
	\rput(3,1.5){\rnode{4}{$\C$}}
	\ncarc{->}{0}{2}\rput(0.9,1.25){$1$}\ncline{->}{0}{1}\Bput*[0.15]{$1$}
	\ncarc{->}{1}{3}\rput(1.2,0.7){$1$}\ncline{->}{1}{2}\Bput*[0.15]{$d$}
	\ncarc{->}{2}{4}\rput(1.8,0.7){$1$}\ncline{->}{2}{3}\Bput*[0.15]{$1$}
	\ncarc{->}{3}{0}\rput(2.05,1.25){$bE$}\ncline{->}{3}{4}\Bput*[0.15]{$d$}
	\ncarc{->}{4}{1}\rput(1.5,1.6){$E$}\ncline{->}{4}{0}\Bput*[0.15]{$dE$}
	}}
	\rput(8.5,1){
	 \scalebox{0.6}{
	\rput(1.5,2.5){\rnode{0}{$\C$}}
	\rput(0,1.5){\rnode{1}{$\C$}}
	\rput(0.5,0){\rnode{2}{$\C$}}
	\rput(2.5,0){\rnode{3}{$\C$}}
	\rput(3,1.5){\rnode{4}{$\C$}}
	\ncarc{->}{0}{2}\rput(0.9,1.25){$1$}\ncline{->}{0}{1}\Bput*[0.15]{$a$}
	\ncarc{->}{1}{3}\rput(1.2,0.7){$1$}\ncline{->}{1}{2}\Bput*[0.15]{$aD$}
	\ncarc{->}{2}{4}\rput(1.8,0.7){$1$}\ncline{->}{2}{3}\Bput*[0.15]{$a$}
	\ncarc{->}{3}{0}\rput(2.1,1.25){$D$}\ncline{->}{3}{4}\Bput*[0.15]{$aD$}
	\ncarc{->}{4}{1}\rput(1.5,1.6){$1$}\ncline{->}{4}{0}\Bput*[0.15]{$aD$}
	}}
\end{pspicture}
\caption{$G$-graphs and corresponding representations of $(Q,R)$ for the group $\frac{1}{5}(1,2)$.}
\label{Z5(1,2)}
\end{center}
\end{figure}
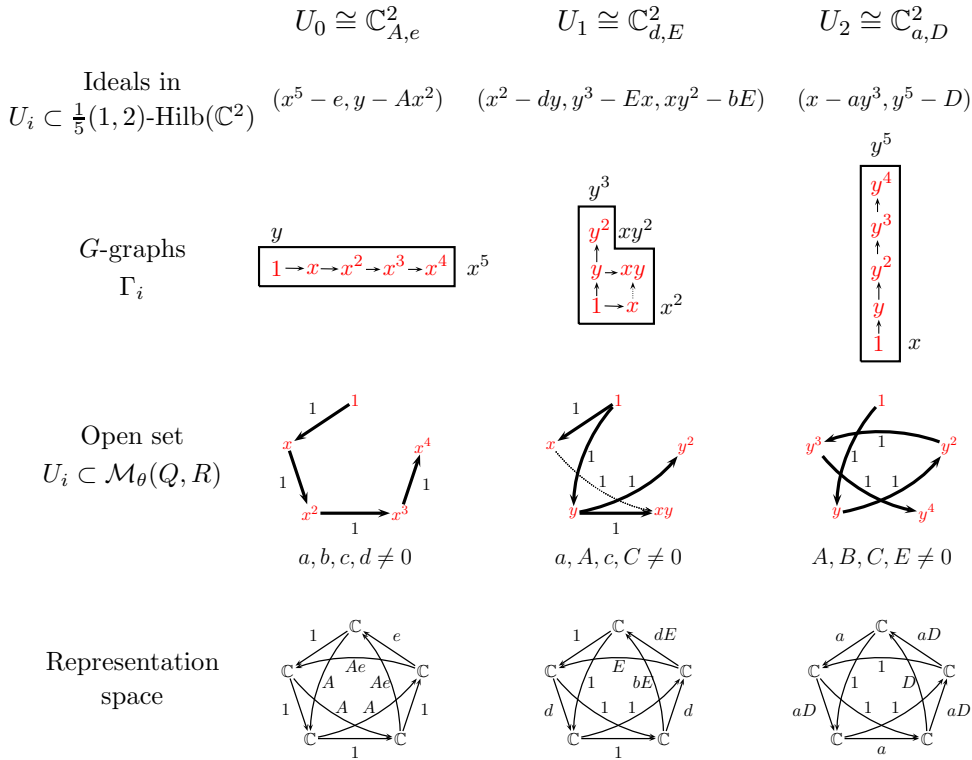

\end{exa}

In the case of $\BD_{2n}(a)$ groups we have the following result.

\begin{thm}\label{opens} Let $G=\BD_{2n}(a)$ and let $\Gamma=\Gamma(r,s;u,v)$ be a $G$-graph with $U_{\Gamma}\subset\mathcal{M}_\theta$. Then,
\begin{itemize}
\item If $\Gamma$ is of type $A$ then $U_{\Gamma_{A}}$ is given by:

$a_{0}, D_{0}, F_{0}\neq0$, and $e_{i}, g_{i}$, $r_{i,j}$, $U'_{i,j}\neq0$ for all $i,j$.

$a_{i}, H_{i}\neq0$ for $i$ even, and $c_{i}, F_{i}\neq0$ for $i$ odd.

For $0<i<u$ set $b_{i},D_{i}\neq0$ if $i$ is even, and $B_{i},d_{i}\neq0$ if $i$ is odd. 

For $i\geq u$ set $B_{i},D_{i}\neq0$.

\item If $\Gamma$ is of type $B$ then $U_{\Gamma_{B}}$ is given by:	

$a_{0}, d_{0}, H_{0}\neq0$, $e_{i}\neq0$, $g_{i}\neq0$, $r_{i,j}\neq0$ for all $i,j$.

$a_{i}, b_{i}, D_{i}, H_{i}\neq0$ for $i$ even, and $B_{i}, c_{i}, d_{i}, F_{i}\neq0$ for $i$ odd. 

If $r>1$ also $C_{0}, R'_{1,1}$, $\ldots$, $R'_{1,r-2}\neq0$.

\begin{itemize}	
\item[(i)] If $\Gamma$ is of type $B_1$ then set $R'_{r+1,1}$,$\ldots$, $R'_{r+1,u-r-2}\neq0$ and $U'_{i,j}\neq0$, $\forall i\neq0,r$ and $\forall j$. For $i=0,r$ we have $U'_{0,r}$,$\ldots$,$U'_{0,q-2}\neq0$ and $U'_{r,u-r},\ldots U'_{r,q-2}\neq0$.

\noindent Also if $q>2$, $C_r\neq0$ if $r$ even, or $A_r\neq0$ if $r$ odd.
	
\item[(ii)] If $\Gamma$ is of type $B_2$ then also set $U'_{i,j}\neq0$, $\forall i>0$ and $\forall j$ and $U'_{0,r},\ldots U'_{0,q-2}\neq0$.
\end{itemize}

\item If $\Gamma$ is of type $C$ then

\begin{itemize}	
\item[(i)] The conditions for $U_{\Gamma_{C^-}}\!\!\subset\mathcal{M}_\theta$ are the same as the ones for $\Gamma_{i}(r,s;q,q)$ with $i=A$ or $B$, and the condition $F_0\neq0$ instead of $H_0\neq0$.  

\item[(ii)] The open conditions for the case $\Gamma_{C^+}$ are the same as those for $\Gamma_{C^-}$ but swapping the conditions for $F_i$ for $H_i$ and vice versa. 
\end{itemize}

\item If $\Gamma$ is of type $D$ then $U_{\Gamma_{D^\pm}}$ is defined by:

\noindent $a_{0}, C_{0}, d_{0}\neq0$, and $a_{i}, b_{i}, D_{i}\neq0$ for $i$ even, $B_{i}, c_{i}, d_{i}\neq0$ for $i$ odd.

\noindent $U'_{i,j}\neq0$ for all $i>0$ and all $j$, $U'_{0,r},\ldots U'_{0,q-2}\neq0$.

\noindent $r_{i,j}\neq0$ for all $i,j$ except for $r_{i,q-i}$, $i=\in[2,k-1]$.

\noindent $R'_{1,1}, \ldots,R'_{1,r-2}\neq0$, $u_{i,q-i}\neq0$ for $i\in[2,k-1]$.

\begin{itemize}	
\item[(i)] If $\Gamma$ is a $G$-graph of type $D^+$ then we also set $e_{0}$, $H_{0}$, $G_{0}$, $E_{1}$, $f_{1}$, $g_{1}$, $H_{1}\neq0$. If $i$ is even then $E_{i}$, $g_{i}$, $F_{i}\neq0$, and if $i$ is odd then $e_{i}$, $G_{i}$, $H_{i}\neq0$.

\item[(ii)] is a $G$-graph of type $D^-$ we set $E_{0}$, $F_{0}$, $g_{0}$, $e_{1}$, $F_{1}$, $G_{1}$, $h_{1}\neq0$. If $i$ is even then $e_{i}$, $G_{i}$, $H_{i}\neq0$, and if $i$ is odd then $E_{i}$, $g_{i}$, $F_{i}\neq0$ with $i\in[0,k-1]$.
\end{itemize}
\end{itemize}
\end{thm}

\begin{proof} An open set in $\mathcal{M}_\theta$ is obtained by making open conditions in the parameter space $\mathbb{V}(I_R)\subset\mathbb{A}^N$. We can change basis at every vertex to take $1$ as basis for every 1-dimensional vertex, and $(1,0)$ and $(0,1)$ for every 2-dimensional. Thus, by \ref{stab-cond} the element $1\in\rho_0^+$ generates the whole representation with this basis. For instance, we always choose $a_0=(1,0)$.

Given any $G$-graph $\Gamma$ the corresponding open set $U_\Gamma\subset\mathcal{M}_\theta$ is obtained by taking the open conditions according to the elements of $\Gamma$. This is done by considering $Q$ to be given (see \ref{section-quiver}) by the $S^G$-modules $S_\rho$ as vertices, and the irreducible maps between them to be the arrows. See Figure \ref{xySegment} for the case $n$ even, where the segment is repeated throughout the quiver. When $n$ is odd replace $e_i$ by $\left(\begin{smallmatrix}x\\-iy\end{smallmatrix}\right)$, $f_i$ by $(y,-ix)$,  $g_i$ by $\left(\begin{smallmatrix}x\\iy\end{smallmatrix}\right)$ and $h_i$ by $(y,ix)$ (to verify the relations in \ref{prop-McKayQ} we have to multiply $r_{i,j}$ and $u_{i,j}$ by $\sqrt2$ for every $i,j$).

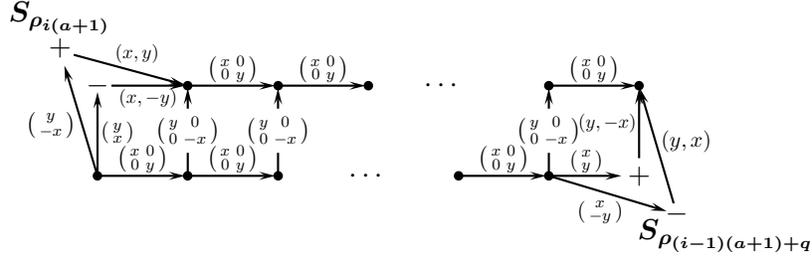
\begin{figure}[h]
\begin{center}
\begin{pspicture}(0,0.4)(10,3.25)
	\psset{nodesep=1pt,arcangle=15,arrowlength=2}
\rput(-1.25,-0.15){
	\rput(2.4,2.4){$-$}\rput(1.9,2.9){+}
	\psdots(3.6,2.4)(4.8,2.4)(6,2.4)(8.4,2.4)(9.6,2.4)
	\rput(7,2.4){$\cdots$}
	\pcline{->}(2.05,2.82)(3.6,2.4)\Aput[0.5pt]{\scalebox{0.65}{$(x,y)$}}
	\pcline{->}(2.55,2.4)(3.6,2.4)\Bput[0.5pt]{\scalebox{0.65}{$(x,-y)$}}
	\pcline{->}(3.6,2.4)(4.8,2.4)
		\Aput[0.5pt]{\scalebox{0.75}{$\left(\begin{smallmatrix}x&0\\0&y\end{smallmatrix}\right)$}}
	\pcline{->}(4.8,2.4)(6,2.4)
		\Aput[0.5pt]{\scalebox{0.75}{$\left(\begin{smallmatrix}x&0\\0&y\end{smallmatrix}\right)$}}
	\pcline{->}(8.4,2.4)(9.6,2.4)
		\Aput[0.5pt]{\scalebox{0.75}{$\left(\begin{smallmatrix}x&0\\0&y\end{smallmatrix}\right)$}}
	\pcline{->}(2.4,1.2)(2.4,2.35)\Bput[0.5pt]{\scalebox{0.75}{$\left(\begin{smallmatrix}y\\x\end{smallmatrix}\right)$}}
	\pcline{->}(2.4,1.2)(1.95,2.7)\Aput[0.1pt]{\scalebox{0.75}{$\left(\begin{smallmatrix}y\\-x\end{smallmatrix}\right)$}}
	\pcline{->}(3.6,1.2)(3.6,2.4)
		\mput*[0.5pt]{\scalebox{0.75}{$\left(\begin{smallmatrix}y&0\\0&-x\end{smallmatrix}\right)$}}
	\pcline{->}(4.8,1.2)(4.8,2.4)
		\mput*[0.5pt]{\scalebox{0.75}{$\left(\begin{smallmatrix}y&0\\0&-x\end{smallmatrix}\right)$}}
	\pcline{->}(8.4,1.2)(8.4,2.4)
		\mput*[0.5pt]{\scalebox{0.75}{$\left(\begin{smallmatrix}y&0\\0&-x\end{smallmatrix}\right)$}}
	\pcline{->}(9.6,1.4)(9.6,2.4)\Aput[0.1pt]{\scalebox{0.65}{$(y,-x)$}}
	\pcline{->}(10.07,0.8)(9.6,2.4)\Bput[0.1pt]{\scalebox{0.75}{$(y,x)$}}

	\psdots(2.4,1.2)(3.6,1.2)(4.8,1.2)(7.2,1.2)(8.4,1.2)
	\rput(9.6,1.2){+}\rput(10.1,0.7){$-$}
	\rput(6,1.2){$\cdots$}
	\pcline{->}(2.4,1.2)(3.6,1.2)
		\Aput[0.5pt]{\scalebox{0.75}{$\left(\begin{smallmatrix}x&0\\0&y\end{smallmatrix}\right)$}}
	\pcline{->}(3.6,1.2)(4.8,1.2)
		\Aput[0.5pt]{\scalebox{0.75}{$\left(\begin{smallmatrix}x&0\\0&y\end{smallmatrix}\right)$}}
	\pcline{->}(7.2,1.2)(8.4,1.2)
		\Aput[0.5pt]{\scalebox{0.75}{$\left(\begin{smallmatrix}x&0\\0&y\end{smallmatrix}\right)$}}
	\pcline{->}(8.4,1.2)(9.4,1.2)\Aput[0.5pt]{\scalebox{0.75}{$\left(\begin{smallmatrix}x\\y\end{smallmatrix}\right)$}}
	\pcline{->}(8.4,1.2)(10,0.72)\Bput[0.1pt]{\scalebox{0.75}{$\left(\begin{smallmatrix}x\\-y\end{smallmatrix}\right)$}}

	\rput(10.75,0.4){$\boldsymbol{S_{\rho_{(i-1)(a+1)+q}}}$}
	\rput(1.9,3.3){$\boldsymbol{S_{\rho_{i(a+1)}}}$}
	}
\end{pspicture}
\caption{Segment $i$ of the quiver between the modules $S_{\rho_{i}}$.}
\label{xySegment}
\end{center}
\end{figure}
By Claim \ref{stab-cond} these irreducible maps send $1\in S_{\rho_{0}^+}$ once to every other $S_{\rho^\pm}$ and twice to every other $S_{V_k}$ linear independently. Denote the polynomials obtained by $f_{\rho^\pm}$ and $(g_V,g_V')$, $(h_V,h_V')$ respectively. In this way, for any stable representation all modules $S_\rho$ have assigned basis polynomials. Thus, if we take the open conditions such that basis elements generated from $1\in S_{\rho_0^+}$ form the $G$-graph $\Gamma$, we obtain the desired open set $U_\Gamma\in\mathcal{M}_\theta$.

If $f\in S_{\rho^\pm}$ and $f\neq f_\rho^\pm$ (i.e.\ $f\notin\Gamma$), then $\exists c\in\C$ such that $f=cf_\rho^\pm$ where $c$ is the path in the representation connecting 1 and $f$ (similarly $(f,f')=c_1(g_V,g_V')+c_2(h_V,h_V')$ for $c_1,c_2\in\C$, $(f,f')\in S_V$). Then $U_\Gamma$ parametrises every $G$-cluster with $\Gamma$ as $G$-graph, so the union of $U_\Gamma$ covers $\mathcal{M}_\theta$. We prove the result case by case. It is worth mentioning that for $\BD_{2n}(a)$ groups we have $(k,q)=1$ (see \cite{thesis} $\S3.3.1$).

\noindent{\em Case A:} We start to generate the representation from $\bf{1}\in\rho_0^+$ and ${\bf a_0}=(1,0)$. We choose to obtain the basis element $(1,0)$ at every 2-dimensional vertex with horizontal arrows taking ${\bf r_{i,j}}=\left(\begin{smallmatrix}1&0\\r'_{i,j}&R'_{i,j}\end{smallmatrix}\right)$ $\forall i,j$, ${\bf a_{i}}=(1,0)$ for $i$ even, ${\bf c_{i}}=(1,0)$ for $i$ odd. The open conditions needed are $r_{i,j}\neq0$ $\forall i,j$, $a_{i}\neq0$ for $i$ even, $c_{i}\neq0$ and $i$ odd. Similarly, we choose to reach $(0,1)$ at every 2-dimensional vertex with vertical arrows taking ${\bf u_{i,j}}=\left(\begin{smallmatrix}u_{i,j}&U_{i,j}\\0&1\end{smallmatrix}\right)$, ${\bf h_{i}}=(0,1)$ for $i$ even, ${\bf f_{i}}=(0,1)$ for $i$ odd, by the open conditions $U'_{i,j}\neq0$ $\forall i,j$, $H_{i}\neq0$ for $i$ even, $F_{i}\neq0$ for $i$ odd (including $i=0$). 

For the 1-dimensional representations on the right hand side we take $e_i=g_i=1$ for all $i$, with the open conditions $e_{i},g_{i}\neq0$. On the left hand side, since the $G$-graph is of type $\Gamma_{A}(r,s;u,v)$ we have $x^{u}y^{u}\notin\Gamma_{A}$ but $x^{i}y^{i}\in\Gamma_{A}$ for $i<u$. In fact, $x^iy^i\in\rho_{i(a+1)}^{(-1)^{i}}$. Thus, we need to reach $\rho_{i(a+1)}^{(-1)^{i}}$ with a nonzero map for $0<i<u$ with a composition of maps of length $i$. We can achieve such a map by taking $d_{1}=b_{2}=d_{3}=1$, $\ldots$ until $d_{u-1}=1$ if $u$ is even, or $b_{u-1}=1$ if $u$ is odd. The condition $x^{u}y^{u}\notin\Gamma_{A}$ is given by $B_{u}=1$ if $u$ is even, or $D_{u}=1$ if $u$ is odd. Finally, from row $u$ to the top row the choices are always $B_{i},D_{i}\neq0$, $i\neq0$ and $D_0\neq0$. 

\noindent {\em Case B:} In this case $x^{u}y^k, x^ky^{u}\notin\Gamma_{B}$, which implies that $x^{i}y^{i}\in\Gamma_{B}$ for $i<u$. This explains the choices at the left hand side of the quiver, while on the right hand side remain the same as before. Since $x^{u}y^m,x^my^{u}\in V_{r}$, the conditions $x^{u}y^m,x^my^{u}\notin\Gamma_{B}$ are expressed with choices $C_{0}, R_{1,1}, R_{1,2},\ldots,R_{1,r-2}\neq0$. If $r\leq k$ we have a $G$-graph of type $B_{1}$, otherwise we have a type $B_{2}$. 

\noindent {\em Case C:} If the $G$-graph $\Gamma(r,s;q,q)$ is of type $B$, then the open conditions are made at the special representation $V_{r}$. The difference between the $C^+$ and $C^-$ is given by $(+)^2\notin\Gamma_{C^+}$ and $(-)^2\notin\Gamma_{C^-}$ which are the choices on the vertical arrows in the right side of $Q$.
 
\noindent {\em Case D:} In this case $(+)\notin\Gamma_{D^+}$ (or $(-)\notin\Gamma_{D^-}$). The open condition is made at the special representation $\rho^+_{q}$ (or at $\rho^-_{q}$ respectively). For instance, in the $D^+$ case we do not allow a path of length $q$ starting from $\rho^+_{0}$ and ending at $\rho^+_{q}$, i.e.\ $E_1=1$. 
\end{proof}

\begin{exa} Let $G=\BD_{12}(7)$ with $q=2$, $k=3$. The $G$-graphs is shown in Figure \ref{GraphsBD12(7)}. 
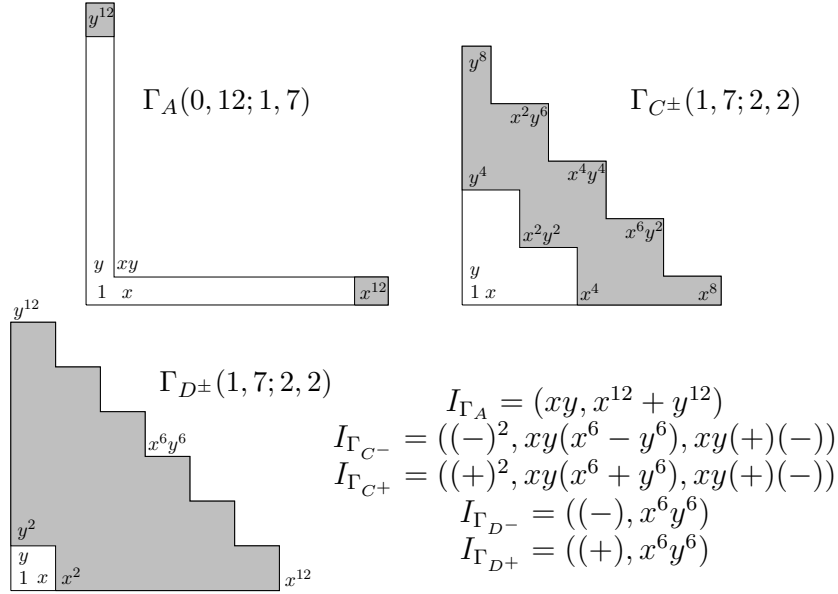
\begin{figure}[ht]
\begin{center}
\begin{pspicture}(0,0.1)(6.5,7.5)
\rput(-3,0){
\rput(2,3.8){
\scalebox{0.85}{
\scalebox{0.35}{
	\psline{-}(0,0)(13.5,0)(13.5,1.25)(1.25,1.25)(1.25,13.5)(0,13.5)(0,0)
	\psline[fillstyle=solid,fillcolor=lightgray](12,0)(13.5,0)(13.5,1.25)(12,1.25)(12,0)
	\psline[fillstyle=solid,fillcolor=lightgray](0,12)(0,13.5)(1.25,13.5)(1.25,12)(0,12)
	}
	\rput(0.25,0.2){\scalebox{0.8}{$1$}}
	\rput(0.65,0.175){\scalebox{0.8}{$x$}}\rput(0.2,0.605){\scalebox{0.8}{$y$}}
	\rput(0.65,0.6){\scalebox{0.8}{$xy$}}
	\rput(4.5,0.225){\scalebox{0.8}{$x^{12}$}}\rput(0.245,4.5){\scalebox{0.8}{$y^{12}$}}
	 }}
\rput(7,3.8){
\scalebox{0.85}{
\scalebox{0.3}{
	\psline{-}(0,0)(13.5,0)(13.5,1.5)(10.5,1.5)(10.5,4.5)(7.5,4.5)(7.5,7.5)(4.5,7.5)(4.5,10.5)(1.5,10.5)(1.5,13.5)(0,13.5)(0,0)
	\psline[fillstyle=solid,fillcolor=lightgray](6,0)(13.5,0)(13.5,1.5)(10.5,1.5)(10.5,4.5)(7.5,4.5)(7.5,7.5)(4.5,7.5)(4.5,10.5)(1.5,10.5)(1.5,13.5)(0,13.5)(0,6)(3,6)(3,3)(6,3)(6,0)
	}
	\rput(0.2,0.2){\scalebox{0.8}{$1$}}
	\rput(0.45,0.17){\scalebox{0.8}{$x$}}\rput(0.2,0.5){\scalebox{0.8}{$y$}}
	\rput(3.85,0.2){\scalebox{0.8}{$x^{8}$}}\rput(0.25,3.8){\scalebox{0.8}{$y^{8}$}}
	\rput(2,0.2){\scalebox{0.8}{$x^4$}}\rput(0.25,2){\scalebox{0.8}{$y^4$}}
	\rput(1.25,1.1){\scalebox{0.8}{$x^2y^2$}}\rput(1.95,2.05){\scalebox{0.8}{$x^4y^4$}}
	\rput(2.85,1.15){\scalebox{0.8}{$x^6y^2$}}\rput(1.05,2.95){\scalebox{0.8}{$x^2y^6$}}
	}}
\rput(1,0){
\scalebox{0.85}{
\scalebox{0.35}{
	\psline{-}(0,0)(12,0)(12,2)(10,2)(10,4)(8,4)(8,6)(6,6)(6,8)(4,8)(4,10)(2,10)(2,12)(0,12)(0,0)
	\psline[fillstyle=solid,fillcolor=lightgray](2,0)(12,0)(12,2)(10,2)(10,4)(8,4)(8,6)(6,6)(6,8)(4,8)(4,10)(2,10)(2,12)(0,12)(0,2)(2,2)(2,0)
	}
	\rput(0.2,0.2){\scalebox{0.8}{$1$}}
	\rput(0.5,0.175){\scalebox{0.8}{$x$}}\rput(0.2,0.5){\scalebox{0.8}{$y$}}
	\rput(0.95,0.225){\scalebox{0.8}{$x^2$}}\rput(0.25,0.95){\scalebox{0.8}{$y^2$}}
	\rput(2.45,2.3){\scalebox{0.8}{$x^6y^6$}}
	\rput(4.5,0.2){\scalebox{0.8}{$x^{12}$}}\rput(0.25,4.4){\scalebox{0.8}{$y^{12}$}}
	 }}
\rput(4,6.5){$\Gamma_{A}(0,12;1,7)$}
\rput(10.5,6.5){$\Gamma_{C^\pm}(1,7;2,2)$}
\rput(4.25,2.75){$\Gamma_{D^\pm}(1,7;2,2)$}
\large{
\rput(8.75,2.5){$I_{\Gamma_{A}}=(xy,x^{12}+y^{12})$}
\rput(8.75,2){$I_{\Gamma_{C^-}}=((-)^2,xy(x^6-y^6),xy(+)(-))$}
\rput(8.75,1.5){$I_{\Gamma_{C^+}}=((+)^2,xy(x^6+y^6),xy(+)(-))$}
\rput(8.75,1){$I_{\Gamma_{D^-}}=((-),x^6y^6)$}
\rput(8.75,0.5){$I_{\Gamma_{D^+}}=((+),x^6y^6)$}
	}}
\end{pspicture}
\caption{The $\BD_{12}(7)$-graphs and the representation ideals, where $(+):=x^2+y^2$ and $(-):=x^2-y^2$.}
\label{GraphsBD12(7)}
\end{center}
\end{figure}
The open choices for $\Gamma_{A}$ and $\Gamma_{C^-}$ are shown in Figure \ref{RepBD12(7)AC-}. Using the quiver in Figure \ref{xySegment} we can calculate the basis polynomials in every irreducible representation (e.g.\ Table \ref{TD12}). The equations of the open cover as hypersurfaces in $\C^3$ are:
\begin{align*}
&U_A: c_0d_1-(c_0d_1^2+1)G_1 &U_{D^+}: e_{1}f_{0}-(e_{1}^2f_{0}-1)D_{1}\\
&U_{C^+}: b_{2}D_{1}-(b_{2}-1)E_{1} &U_{D^-}: g_{1}h_{0}-(g_{1}^2h_{0}-1)D_{1}\\
&U_{C^-}: b_2D_1-(b_{2}-1)G_{1}\\
\end{align*} 
The dual graph of the exceptional divisor in $\Hilb{G}$ is $\begin{smallmatrix}-2\\-2\end{smallmatrix}\text{\small{$-3$\ $-2$}}$.

{\renewcommand{\arraystretch}{1}
\begin{table}[h]
{\small 
\begin{center}
\begin{tabular}{|c|l|c|l|c|l|}
\hline
\!\!$\rho_{0}^+$\!\! & 1			& \!\!$\rho_{6}^+$\!\!  & $(+)^3$  	& \!\!$V_1$\!\!	& $(1,0)=(x,y)$ \\
\!\!$\rho_{0}^-$\!\!  & $2xy(+)^2$ 	& \!\!$\rho_{6}^-$\!\!  & $2xy(+)^5$	& 		& $(0,1)=(y(+)^3, x(+)^3)$ \\
\!\!$\rho_{2}^+$\!\! & $(+)$			& \!\!$\rho_{8}^+$\!\!  & $(+)^4$  	& \!\!$V_5$\!\!	& $(1,0)=(x(+)^2,y(+)^2)$ \\	
\!\!$\rho_{2}^-$\!\!  & $2xy(+)^3$	& \!\!$\rho_{8}^-$\!\!  & $2xy$ 		& 		& $(0,1)=(y(+)^5,-x(+)^5)$ \\
\!\!$\rho_{4}^+$\!\!  & $(+)^2$ 		& \!\!$\rho_{10}^+$\!\!  & $(+)^5$ 	& \!\!$V_9$\!\! 	& $(1,0)=(x(+)^4,y(+)^4)$ \\
\!\!$\rho_{4}^-$\!\!  & $2xy(+)^4$ 	& \!\!$\rho_{10}^-$\!\!	& $2xy(+)$	& 		& $(0,1)=(y(+),-x(+))$ \\
\hline
\end{tabular}
\caption{Basis elements of the $G$-graph $\Gamma_{D^-}(1,7;2,2)$.}
\label{TD12}
\end{center}}
\end{table}}

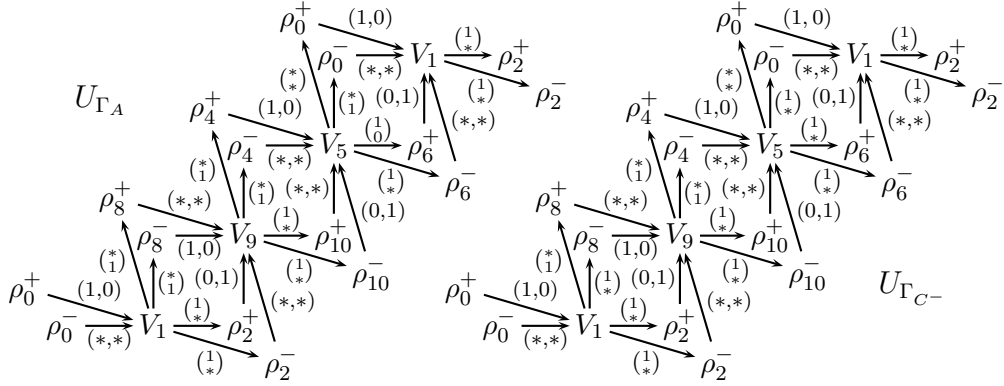
\begin{figure}[ht]
\begin{center}
\begin{pspicture}(0.5,-1.85)(5.75,2.1)
\psset{nodesep=2pt}
\rput(-3.8,-1.25){
	\rput(0.8,0){
	\rput(-0.5,0.5){\rnode{0+}{$\rho_0^+$}}
	\rput(0,0){\rnode{0-}{$\rho_0^-$}}
	\rput(1.2,0){\rnode{1}{$V_1$}}
	\rput(2.4,0){\rnode{2+}{$\rho_2^+$}}
	\rput(2.9,-0.5){\rnode{2-}{$\rho_2^-$}}
	\ncline{->}{0+}{1}\Aput[0.1pt]{\scriptsize $(1,\!0)$}
	\ncline{->}{0-}{1}\Bput[0.5pt]{\scriptsize $(\ast,\!\ast)$}
	\ncline{->}{1}{2+}\Aput[0.5pt]{\footnotesize $\left(\!\begin{smallmatrix}1\\\ast\end{smallmatrix}\!\right)$}
	\ncline{->}{1}{2-}\Bput[0.5pt]{\footnotesize $\left(\!\begin{smallmatrix}1\\\ast\end{smallmatrix}\!\right)$}
	\rput(0.7,1.7){\rnode{8+}{$\rho_8^+$}}
	\rput(1.2,1.2){\rnode{8-}{$\rho_8^-$}}
	\rput(2.4,1.2){\rnode{9}{$V_9$}}
	\rput(3.6,1.2){\rnode{10+}{$\rho_{10}^+$}}
	\rput(4.1,0.7){\rnode{10-}{$\rho_{10}^-$}}
	\ncline{->}{8+}{9}\Aput[0pt]{\scriptsize $(\ast,\!\ast)$}
	\ncline{->}{8-}{9}\Bput[0.5pt]{\scriptsize $(1,\!0)$}
	\ncline{->}{9}{10+}\Aput[0.5pt]{\footnotesize $\left(\!\begin{smallmatrix}1\\\ast\end{smallmatrix}\!\right)$}
	\ncline{->}{9}{10-}\Bput[0.1pt]{\footnotesize $\left(\!\begin{smallmatrix}1\\\ast\end{smallmatrix}\!\right)$}
	\rput(1.9,2.9){\rnode{4+}{$\rho_4^+$}}
	\rput(2.4,2.4){\rnode{4-}{$\rho_4^-$}}
	\rput(3.6,2.4){\rnode{5}{$V_5$}}
	\rput(4.8,2.4){\rnode{6+}{$\rho_{6}^+$}}
	\rput(5.3,1.9){\rnode{6-}{$\rho_{6}^-$}}
	\ncline{->}{4+}{5}\Aput[0.5pt]{\scriptsize $(1,\!0)$}
	\ncline{->}{4-}{5}\Bput[0.5pt]{\scriptsize $(\ast,\!\ast)$}
	\ncline{->}{5}{6+}\Aput[0.5pt]{\footnotesize $\left(\!\begin{smallmatrix}1\\0\end{smallmatrix}\!\right)$}
	\ncline{->}{5}{6-}\Bput[0.1pt]{\footnotesize $\left(\!\begin{smallmatrix}1\\\ast\end{smallmatrix}\!\right)$}
	\rput(3.1,4.1){\rnode{00+}{$\rho_0^+$}}
	\rput(3.6,3.6){\rnode{00-}{$\rho_0^-$}}
	\rput(4.8,3.6){\rnode{11}{$V_1$}}
	\rput(6,3.6){\rnode{22+}{$\rho_{2}^+$}}
	\rput(6.5,3.1){\rnode{22-}{$\rho_{2}^-$}}
	\ncline{->}{00+}{11}\Aput[0.5pt]{\scriptsize $(1,\!0)$}
	\ncline{->}{00-}{11}\Bput[0.5pt]{\scriptsize $(\ast,\!\ast)$}
	\ncline{->}{11}{22+}\Aput[0.5pt]{\footnotesize $\left(\!\begin{smallmatrix}1\\\ast\end{smallmatrix}\!\right)$}
	\ncline{->}{11}{22-}\Bput[0.1pt]{\footnotesize $\left(\!\begin{smallmatrix}1\\\ast\end{smallmatrix}\!\right)$}

	\ncline{->}{1}{8-}\Bput[0.5pt]{\footnotesize $\left(\!\begin{smallmatrix}\ast\\1\end{smallmatrix}\!\right)$}
	\ncline{->}{9}{4-}\Bput[0.5pt]{\footnotesize $\left(\!\begin{smallmatrix}\ast\\1\end{smallmatrix}\!\right)$}
	\ncline{->}{5}{00-}\Bput[0.5pt]{\footnotesize $\left(\!\begin{smallmatrix}\ast\\1\end{smallmatrix}\!\right)$}

	\ncline{->}{1}{8+}\aput[0.5pt](0.6){\footnotesize $\left(\!\begin{smallmatrix}\ast\\1\end{smallmatrix}\!\right)$}
	\ncline{->}{9}{4+}\aput[0.5pt](0.6){\footnotesize $\left(\!\begin{smallmatrix}\ast\\1\end{smallmatrix}\!\right)$}
	\ncline{->}{5}{00+}\aput[0.5pt](0.6){\footnotesize $\left(\!\begin{smallmatrix}\ast\\\ast\end{smallmatrix}\!\right)$}

	\ncline{->}{2+}{9}\Aput[0.5pt]{\scriptsize $(0,\!1)$}
	\ncline{->}{10+}{5}\Aput[0.5pt]{\scriptsize $(\ast,\!\ast)$}
	\ncline{->}{6+}{11}\Aput[0.5pt]{\scriptsize $(0,\!1)$}

	\ncline{->}{2-}{9}\bput[0.5pt](0.4){\scriptsize $(\ast,\!\ast)$}
	\ncline{->}{10-}{5}\bput[0.5pt](0.4){\scriptsize $(0,\!1)$}
	\ncline{->}{6-}{11}\bput[0.5pt](0.4){\scriptsize $(\ast,\!\ast)$}

	\rput(0.5,3){\large $U_{\Gamma_{A}}$}
	}}
	
	\rput(2,-1.25){
	\rput(0.8,0){
	\rput(-0.5,0.5){\rnode{0+}{$\rho_0^+$}}
	\rput(0,0){\rnode{0-}{$\rho_0^-$}}
	\rput(1.2,0){\rnode{1}{$V_1$}}
	\rput(2.4,0){\rnode{2+}{$\rho_2^+$}}
	\rput(2.9,-0.5){\rnode{2-}{$\rho_2^-$}}
	\ncline{->}{0+}{1}\Aput[0.1pt]{\scriptsize $(1,\!0)$}
	\ncline{->}{0-}{1}\Bput[0.5pt]{\scriptsize $(\ast,\!\ast)$}
	\ncline{->}{1}{2+}\Aput[0.5pt]{\footnotesize $\left(\!\begin{smallmatrix}1\\\ast\end{smallmatrix}\!\right)$}
	\ncline{->}{1}{2-}\Bput[0.5pt]{\footnotesize $\left(\!\begin{smallmatrix}1\\\ast\end{smallmatrix}\!\right)$}
	\rput(0.7,1.7){\rnode{8+}{$\rho_8^+$}}
	\rput(1.2,1.2){\rnode{8-}{$\rho_8^-$}}
	\rput(2.4,1.2){\rnode{9}{$V_9$}}
	\rput(3.6,1.2){\rnode{10+}{$\rho_{10}^+$}}
	\rput(4.1,0.7){\rnode{10-}{$\rho_{10}^-$}}
	\ncline{->}{8+}{9}\Aput[0pt]{\scriptsize $(\ast,\!\ast)$}
	\ncline{->}{8-}{9}\Bput[0.5pt]{\scriptsize $(1,\!0)$}
	\ncline{->}{9}{10+}\Aput[0.5pt]{\footnotesize $\left(\!\begin{smallmatrix}1\\\ast\end{smallmatrix}\!\right)$}
	\ncline{->}{9}{10-}\Bput[0.1pt]{\footnotesize $\left(\!\begin{smallmatrix}1\\\ast\end{smallmatrix}\!\right)$}
	\rput(1.9,2.9){\rnode{4+}{$\rho_4^+$}}
	\rput(2.4,2.4){\rnode{4-}{$\rho_4^-$}}
	\rput(3.6,2.4){\rnode{5}{$V_5$}}
	\rput(4.8,2.4){\rnode{6+}{$\rho_{6}^+$}}
	\rput(5.3,1.9){\rnode{6-}{$\rho_{6}^-$}}
	\ncline{->}{4+}{5}\Aput[0.5pt]{\scriptsize $(1,\!0)$}
	\ncline{->}{4-}{5}\Bput[0.5pt]{\scriptsize $(\ast,\!\ast)$}
	\ncline{->}{5}{6+}\Aput[0.5pt]{\footnotesize $\left(\!\begin{smallmatrix}1\\\ast\end{smallmatrix}\!\right)$}
	\ncline{->}{5}{6-}\Bput[0.1pt]{\footnotesize $\left(\!\begin{smallmatrix}1\\\ast\end{smallmatrix}\!\right)$}
	\rput(3.1,4.1){\rnode{00+}{$\rho_0^+$}}
	\rput(3.6,3.6){\rnode{00-}{$\rho_0^-$}}
	\rput(4.8,3.6){\rnode{11}{$V_1$}}
	\rput(6,3.6){\rnode{22+}{$\rho_{2}^+$}}
	\rput(6.5,3.1){\rnode{22-}{$\rho_{2}^-$}}
	\ncline{->}{00+}{11}\Aput[0.5pt]{\scriptsize $(1,0)$}
	\ncline{->}{00-}{11}\Bput[0.5pt]{\scriptsize $(\ast,\!\ast)$}
	\ncline{->}{11}{22+}\Aput[0.5pt]{\footnotesize $\left(\!\begin{smallmatrix}1\\\ast\end{smallmatrix}\!\right)$}
	\ncline{->}{11}{22-}\Bput[0.1pt]{\footnotesize $\left(\!\begin{smallmatrix}1\\\ast\end{smallmatrix}\!\right)$}

	\ncline{->}{1}{8-}\Bput[0.5pt]{\footnotesize $\left(\!\begin{smallmatrix}1\\\ast\end{smallmatrix}\!\right)$}
	\ncline{->}{9}{4-}\Bput[0.5pt]{\footnotesize $\left(\!\begin{smallmatrix}\ast\\1\end{smallmatrix}\!\right)$}
	\ncline{->}{5}{00-}\Bput[0.5pt]{\footnotesize $\left(\!\begin{smallmatrix}1\\\ast\end{smallmatrix}\!\right)$}

	\ncline{->}{1}{8+}\aput[0.5pt](0.6){\footnotesize $\left(\!\begin{smallmatrix}\ast\\1\end{smallmatrix}\!\right)$}
	\ncline{->}{9}{4+}\aput[0.5pt](0.6){\footnotesize $\left(\!\begin{smallmatrix}\ast\\1\end{smallmatrix}\!\right)$}
	\ncline{->}{5}{00+}\aput[0.5pt](0.6){\footnotesize $\left(\!\begin{smallmatrix}\ast\\\ast\end{smallmatrix}\!\right)$}

	\ncline{->}{2+}{9}\Aput[0.5pt]{\scriptsize $(0,\!1)$}
	\ncline{->}{10+}{5}\Aput[0.5pt]{\scriptsize $(\ast,\!\ast)$}
	\ncline{->}{6+}{11}\Aput[0.5pt]{\scriptsize $(0,\!1)$}

	\ncline{->}{2-}{9}\bput[0.5pt](0.4){\scriptsize $(\ast,\!\ast)$}
	\ncline{->}{10-}{5}\bput[0.5pt](0.4){\scriptsize $(0,\!1)$}
	\ncline{->}{6-}{11}\bput[0.5pt](0.4){\scriptsize $(\ast,\!\ast)$}
	\rput(5.5,0.5){\large $U_{\Gamma_{C^-}}$}
	}}
\end{pspicture}
\caption{Open sets $U_{\Gamma_{A}}$, $U_{\Gamma_{C^-}}\subset\mathcal{M}_\theta$ for $\BD_{12}(7)$.}
\label{RepBD12(7)AC-}
\end{center}
\end{figure}

\end{exa}


\end{document}